\documentclass[twoside,11pt]{amsart}
\usepackage{amsmath,latexsym,amssymb,times,enumerate}

\usepackage{verbatim}

\def\demo{\noindent{\bf Proof. }}
\def\sqr#1#2{{\vcenter{\hrule height.#2pt
        \hbox{\vrule width.#2pt height#1pt \kern#1pt
                \vrule width.#2pt}
        \hrule height.#2pt}}}
\def\square{\mathchoice\sqr64\sqr64\sqr{4}3\sqr{3}3}
\def\QED{\hfill$\square$}

\def\lto{\longrightarrow}
\def\CC{\mathbb C}
\def\Z{\mathbb Z}
\def\P{\mathbb P}
\def\Q{\mathbb Q}

\def\oa{{\overline a}}
\def\ob{{\overline b}}
\def\o1{{\overline 1}}
\def\m{{\mathfrak m}}

\def\a{{\mathfrak a}}
\def\b{{\mathfrak b}}
\def\c{{\mathfrak c}}
\def\e{{\mathfrak e}}

\def\r{\rho}
\def\ms{\medskip}
\def\bs{\bigskip}
\def\s{\smallskip}

\def\C{\mathcal C}
\def\E{\mathcal E}
\def\L{\mathcal L}
\def\N{\mathcal N}
\def\O{\mathcal O}
\def\Pic{{\rm Pic}\, }
\def\Apic{{\rm A Pic} \,}
\def\Spec{{\rm Spec} \,}
\def\Cart{{\rm Cart} \,}
\def\I{\mathcal I}
\def\nC{\widetilde{\mathcal C}}

\def \d{\mathfrak d}
\def \cA{\widehat{A}}
\def \cma{\widehat{\mathfrak a}}
\def\nA{\widetilde{A}}
\def\nD{\widetilde{D}}
\def\G{\Bbb G}
\def\Div{\rm Div \,}

\newtheorem{Theorem}{Theorem}[section]
\newtheorem{Lemma}[Theorem]{Lemma}
\newtheorem{Corollary}[Theorem]{Corollary}
\newtheorem{Proposition}[Theorem]{Proposition}

\newtheorem{Notation and Discussion}[Theorem]{Notation and Discussion}
\newtheorem{Assumptions and Discussion}[Theorem]{Assumptions and Discussion}
\newtheorem{Assumptions}[Theorem]{Assumptions}
\newtheorem{Remark}[Theorem]{Remark}
\newtheorem{Example}[Theorem]{Example}
\newtheorem{Examples}[Theorem]{Examples}
\newtheorem{Note}[Theorem]{Note}
\newtheorem{Definition}[Theorem]{Definition}

\begin{document}

\baselineskip=16pt

\title[Divisors class groups of singular surfaces]
{\Large\bf Divisors class groups of singular surfaces}

\author[R. Hartshorne \and C. Polini]
{Robin Hartshorne  \and Claudia Polini}

\thanks{AMS 2010 {\em Mathematics Subject Classification}.
Primary 14C20, 13A30; Secondary 14M10, 14J05.}

\thanks{The second author was partially supported by the
NSA and the NSF}

\address{Department of Mathematics, University of California, Berkeley, California, 94720-3840}\email{robin@math.berkeley.edu}
\address{Department of Mathematics, University of Notre Dame, Notre Dame, Indiana 46556} \email{cpolini@nd.edu}

\vspace{-0.1in}

\begin{abstract} We compute divisors class groups of singular surfaces. Most notably we produce an exact sequence that relates the Cartier divisors and almost Cartier divisors of a surface to the those of its normalization. This generalizes Hartshorne's theorem for the cubic ruled surface in $\mathbb P^3$. We apply these results to limit the possible curves that can be set-theoretic complete intersection in $\P^3$ in characteristic zero. \end{abstract}

\maketitle

\vspace{-0.2in}

\section{Introduction}



On a nonsingular variety, the study of divisors and linear systems is classical. In fact the entire theory of curves and surfaces is dependent on this study of codimension one subvarieties and the linear and algebraic families in which they move. 

This theory has been generalized in two directions:  the Weil divisors on a normal variety, taking codimension one subvarieties as prime divisors; and the Cartier divisors on an arbitrary scheme, based on locally principal codimension one subschemes. Most of the literature both in algebraic geometry and commutative algebra up to now has been limited to these kinds of divisors. 

More recently there have been good reasons to consider divisors on non-normal varieties. Jaffe \cite{Jaffe} introduced the notion of an almost Cartier divisor, which is locally principal off a subset of codimension two. A theory of generalized divisors was proposed on curves in \cite{H0}, and extended to any dimension in \cite{H}. The latter paper gave a complete description of the generalized divisors on the ruled cubic surface in $\P^3$.

In this paper we extend that analysis to an arbitrary integral surface $X$, explaining the group $\Apic X$ of linear equivalence classes of almost Cartier divisors on $X$ in terms of the Picard group of the normalization $S$ of $X$ and certain local data at the singular points of $X$. We apply these results to give limitations on the possible curves that can be set-theoretic compete intersections in $\P^3$ in characteristic zero 

In section 2 we explain our basic set-up, comparing divisors on a variety $X$ to its normalization $S$. In Section 3 we prove a local isomorphism that computes the group of almost Cartier divisors at a singular point of  $X$ in terms of the Cartier divisors along the curve of singularities and its inverse image in the normalization. In Section 4 we derive some global exact sequences for the groups $\Pic X$, $\Apic X$, and $\Pic S$, which generalize the results of \cite[\S 6]{H} to arbitrary surfaces These results are particularly transparent for surfaces with ordinary singularities, meaning a double curve with a finite number of pinch points and triple points. 

In section 5 we gather some results on curves that we need in our calculations on surfaces. Then in Section 6 we give a number of examples of surfaces and compute their groups of almost Picard divisors. 

In Section 7 we apply these results to limit the possible degree and genus of curves in $\P^3$ that can be set-theoretic complete intersections on surfaces with ordinary singularities in characteristic zero, extending earlier work of Jaffe and Boratynski. We illustrate these results with the determination of all set-theoretic complete intersections on a number of particular surfaces in $\P^3$. 

Our main results assume that the ground field $k$ is of characteristic zero, so that a) we can use the exponential sequence in comparing the additive and multiplicative structures, and b) so that the additive group of the field is a torsion-free abelian group.

The first author would like to thank the Department of Mathematics at the University of Notre Dame for hospitality during the preparation of this paper.

\section{Divisors and Finite Morphisms}
\s

All the rings treated in this paper are Noetherian,  essentially of finite type over a field $k$ which is algebraically closed.  
In our application we will often compare divisors on integral surface $X$ with its normalization $S$. But some of our preliminary results are valid more generally so we fix a set of assumptions. 
\begin{Assumptions}\label{A1} Let $\pi:S \longrightarrow X$ be a dominant finite morphism of reduced schemes. Let $\Gamma$ and $L$ be codimension one subschemes in $S$ and $X$ respectively 
such that $\pi$ restricts to a morphism of $\Gamma$ to $L$.
Assume that $S$ and $X$ both satisfy $G_1$ (i.e. Gorenstein in codimension 1) and $S_2$ (i.e. Serre's condition $S_2$) so that the theory of generalized divisors developed in \cite{H} can be applied. Further assume that the schemes $\Gamma$ and $L$ have no embedded associated primes, hence they satisfy $S_1$.
 \end{Assumptions}

Now we recall the notion of generalized divisors from \cite{H}. If $X$ is a scheme satisfying $G_1$ and $S_2$, we denote by $K_X$ the sheaf of total quotient rings of the structure sheaf ${\mathcal O}_X$.  A {\it generalized divisor} on $X$ is a {\it fractional ideal} ${\mathcal I} \subset K_X$, i.e.  a coherent sub-$\mathcal O_X$-module of $K_X$, that is {\it nondegenerate}, namely for each generic point $\eta \in X$, $\mathcal I_{\eta} =K_{X,\eta}$, and such  that $\mathcal I$  is a reflexive $\mathcal O_X$-module. 

We say $\mathcal I$ is {\it principal} if it is generated by a single non-zero-divisor $f$ in $K_X$. We say $\mathcal I$ is {\it Cartier} if it is locally 
principal everywhere. We say $\mathcal I$ is {\it almost Cartier} if it is locally principal off subsets of codimension at least 2.  We denote by ${\rm Cart} X$ and by ${\rm A Cart} X$
the groups of Cartier divisors and almost Cartier divisors, respectively, and dividing these by the subgroup of principal divisors we obtain the divisors class groups ${\rm Pic} X$ and ${\rm A Pic} X$, respectively. The divisor $\mathcal I$ is {\it effective} if it is contained in $\mathcal O_X$. In that case it defines a codimension one subscheme $Y \subset X$ without embedded components. Conversely, for any such $Y$, its sheaf of ideals $\mathcal I_{Y}$ is an effective divisor. 

We recall some properties of these groups.

\begin{Proposition}\label{P1} Adopt assumptions \ref{A1}. The following hold:
\begin{enumerate}
\item[(a)] There is a natural map $\pi^{\star}: {\rm Pic} \, X \longrightarrow {\rm Pic} \, S$
\item[(b)] There is a natural map $\pi^{\star}: {\rm A Pic} \, X \longrightarrow {\rm A Pic} \, S$
\item[(c)] There is an exact sequence
$$ 0 \rightarrow {\rm Pic} X \longrightarrow {\rm A Pic} \, X \longrightarrow \bigoplus_{x \in X} {\rm A Pic} ({\rm Spec} \, \mathcal O_{X,x})
$$
where the sum is taken over all points $x \in X$ of codimension at least two.
\end{enumerate}
\end{Proposition}
\demo For $(a)$ and $(b)$ see \cite[2.18]{H}. The map on ${\rm Pic}$ makes sense for any morphism of schemes. For ${\rm A Pic}$, we need only to observe that since $\pi$ is a dominant finite morphism, if $Z\subset X$ has codimension two, then also $\pi^{-1}(Z) \subset S$ has codimension two.   
The sequence in  $(c)$ is due to Jaffe for surfaces (see \cite[2.15]{H}), but holds in any dimension (same proof). 
\QED

\begin{Proposition}\label{P2} Adopt assumptions \ref{A1}. Further assume that $X$ and $S$ are affine and $S$ is smooth.  Then there is a natural group homomorphism 
$$ \varphi: {\rm A Pic} X \longrightarrow {\rm Cart} \, \Gamma/ \pi^{*}{\rm Cart} \, L
$$
\end{Proposition}
\demo Given a divisor class $\d \in {\rm A Pic} \, X$, choose an effective divisor $D \in \d$ that does not contain any irreducible component of $L$ in its support (this is possible by Lemma \ref{2.4} below). Now restrict the divisor $D$ to $X-L$, transport it via the isomorphism $\pi$ to $S-\Gamma$, and take its closure in $S$. Since $S$ is smooth, this will be a Cartier divisor on all of $S$, which we can intersect with $\Gamma$ to give a Cartier divisor on $\Gamma$. 

If we choose another effective divisor $D'$ representing the same class $\d$, that also does not contain any component of $L$ in its support, then $D-D'$ is a principal divisor $(f)$ for some $f \in K_X$.  Since $\pi$ gives an isomorphism of $S-\Gamma$ to $X-L$ it follows that $S$ and $X$ are birational, i.e.  $K_X=K_S$. So the equation $D'-D=(f)$ persists on $S$, showing that the ambiguity of our construction is the Cartier divisor on $\Gamma$ defined by the restriction of $f$. Note now that since $(f)=D-D'$, we can write $\mathcal I_{D'}=f\mathcal I_D$ where $\mathcal I_{D'}$ and $\mathcal I_D$ are the ideals of $D'$ and $D$ in $\mathcal O_X$, and $f\in K_X$. If $\lambda$ is a generic point of $L$, then after localizing, the ideals ${\mathcal I}_{D',\lambda}$ and ${\mathcal I}_{D,\lambda}$ are both the whole ring $\mathcal O_{X,\lambda}$, since $D$ and $D'$ are effective divisors not containing any component of $L$ in their support. Therefore $f$ is a unit in $O_{X,\lambda}$. Thus $f$ restricts to a non-zerodivisor in the total quotient ring $K_L$, whose stalk at $\lambda$ is isomorphic to $O_{X,\lambda}/{\mathcal I}_{L,\lambda}$. Thus the restriction of $f$ defines a Cartier divisor on $L$ whose image in $\Gamma$ will be the same as the restriction of $f$ from $S$ to $\Gamma$. Hence our map $\varphi$ is well-defined to the quotient group $ {\rm Cart} \, \Gamma/ \pi^{*}{\rm Cart} \, L$.
\QED

\bs

The following lemma is the affine analogue of \cite[2.11]{H}.

\s

\begin{Lemma}\label{2.4} Let $X$ be an affine scheme satisfying $G_1$ and $S_2$. Let ${\mathfrak d} \in {\rm A Pic} \, X$ be an equivalence class of almost Cartier divisors. Let $Y_1, \ldots, Y_r$ be irreducible codimension one subsets of $X$. Then there exists an effective divisor $D \in \d$ that contains none of the $Y_i$ in its support. 
\end{Lemma} 
\demo The class $\d$ corresponds to a reflexive coherent sheaf $\mathcal L$ of $\mathcal O_X$-modules \cite[2.8]{H}, which is locally free at all points $x \in X$ of codimension one because the divisors in $\d$ are almost Cartier. Since $X$ is affine, the sheaf $\mathcal L$ is generated by global sections. Thus for the generic point $y_i$ of $Y_i$, there will be a section $s_i \in \Gamma(X, \mathcal L)$ whose image in the stalk $\mathcal L_{y_i}$ is not contained in $\m_{y_i}\mathcal L_{y_i}$. Those sections $s \in \Gamma(X,\mathcal L)$ not having this property form a proper sub-vector space $V_i$ of $\Gamma(X,\mathcal L)$. Now if we choose a section $s \in \Gamma(X,\mathcal L)$ not contained in any of the $V_i$, the corresponding divisor $D$ \cite[2.9]{H} will be an effective divisor in the class $\d$, not containing any of the $Y_i$ in its support.
\QED

\bs
\begin{Remark}\label{2.5}{\rm  In Proposition \ref{P2} if $S$ is not smooth  the construction does not work because the closure of $D$ in $S$ may not be Cartier. However if we define $G$ to be the following subset of ${\rm A Pic X}$, namely
$$G = \{  \d \in {\rm A Pic} \, X \mid \pi^{*}(\d) \in {\rm Pic} \, S \}$$
then we can construct the map $\varphi: G \longrightarrow  {\rm Cart} \, \Gamma/ \pi^{*}{\rm Cart} \, L$ in the same way. The condition that the element $\pi^{*}(\d)$ of ${\rm A Pic} \, S$ lies in $ {\rm Pic} \, S$ is equivalent, by Proposition \ref{P1}(c), to the vanishing of its image in $ {\rm A Pic} ({\rm Spec} \, \mathcal O_{S,s})$ for all singular points $s \in S$.}
 
\end{Remark}

\bs

\begin{Proposition}\label{2.6} 
Adopt assumptions \ref{A1}. Assume that the map induced by $\pi$ from $\mathcal I_{L,X}$ to $  \pi_{*} (\mathcal I_{\Gamma, X})$ is an isomorphism. Then the map of sheaves of abelian groups $\gamma: {\mathcal N} \longrightarrow {\mathcal N_0}$ on $X$ defined by the following diagram is an isomorphism:
\begin{equation}\label{eq1}
\begin{array}{ccccccc}
 \mathcal O_{X}^{*}  & \longrightarrow &  \pi_*\mathcal O_{S}^{*}  & \longrightarrow
& \mathcal N & \rightarrow & 0 \\
\downarrow\scriptstyle{\alpha} &  & \downarrow\scriptstyle{\beta}  &  & \downarrow\scriptstyle{\gamma}  &  &  \\
 \mathcal O_{L}^{*} & \longrightarrow & \pi_*\mathcal O_{\Gamma}^{*} & \longrightarrow
& \mathcal N_0 & \rightarrow & 0 \\
\end{array}
\end{equation}
\end{Proposition}
\demo  For every point $x \in X$,  set $(A, \m_A)$ to be the local ring $\mathcal O_{X,x}$ and $B$ to be the semi-local ring $\mathcal O_{S, \pi^{-1}(x)}$.  As it is sufficient to check an isomorphism of sheaves on stalks,  we can restrict to the local situation where $X={\rm Spec} \, A$ and $S={\rm Spec} \, B$. Let $A_0=A/I$ be the local ring of $L$ and $B_0=B/J$ be the semi-local ring of $\Gamma$. Our hypothesis says that the homomorphism  from $A$ to $B$  induces an isomorphism from  $I$ to $J$. Now we consider the diagram of abelian groups:
\begin{equation*}
\begin{array}{ccccccc}
 A^{*}  & \longrightarrow &  B^{*}  & \longrightarrow
& N & \rightarrow & 0 \\
 \downarrow\scriptstyle{\a} &  & \downarrow\scriptstyle{\b}  &  & \downarrow\scriptstyle{\c}&  &  \\
 A_0^{*} & \longrightarrow & B_0^{*} & \longrightarrow
&  N_0 & \rightarrow & 0 \\
\end{array}
\end{equation*}
and we need to show that the induced map $\c$ is an isomorphism.

Since $A \longrightarrow A_0$ and $B \longrightarrow B_0$ are surjective maps of (semi)-local rings, the corresponding maps on units  $\a$ and $\b$ are surjective (see Lemma  \ref{units}). Therefore the third map $\c$ is surjective.

To show $\c$ is injective, let $a \in N$ go to $1$ in $N_0$. Because the diagram is commutative $a$ comes from an element  $b \in B^{*}$ and $\b(b)=c \in   B_0^{*}$ whose image in $N_0$ is $1$. Hence $c$ comes from an element  $d\in A_0^{*} $, which lifts to an $e\in A^{*}$. Let  $f $ be the image of $e$ in $B^{*}$. Now $b$ and $f$ have the same image $c$ in $B^{*}_0$. Regarding them as elements of the ring $B$ this means that their difference is in the ideal $J$. But $J$ by hypothesis is isomorphic to $I$, hence there is an element $g \in I$ whose image gives $b-f$ in $J$. Now consider the element $h=g+e \in A$. Since $g \in I \subset \m_A$, the element $h$ is a unit in $A$, i.e. it is an element of $A^{*}$. Furthermore its image in $B^{*}$ is $b$. Therefore the image of $b$ in $N$, which is  $a$, is equal to $1$. Thus the map $\c$ is an isomorphism. Since this holds at all stalks $x\in X$, we conclude that the map $\gamma: \mathcal N \longrightarrow \mathcal N_0$ of sheaves is an isomorphism. 
\QED

\bs

\begin{Remark}\label{R2}{\rm  In applications we will often consider a situation where $X$ is an integral scheme and $S$ is its normalization. Then ${\mathcal O}_S$ is a generalized divisor on $X$, whose inverse $\mathcal I=\{ a \in {\mathcal O_X} \mid a {\mathcal O_S} \subset { \mathcal O_X} \}$ is just the {\it conductor} of the integral extension. If we define $L$ by this ideal, and $\Gamma$ as $\pi^{*}(L)$, then   the map induced by $\pi$ from $\mathcal I_{L,X}$ to $  \pi_{*} (\mathcal I_{\Gamma, X})$ is an isomorphism and the hypothesis on the ideal sheaves is satisfied. Conversely, if the map from $\mathcal I_{L,X}$ to $  \pi_{*} (\mathcal I_{\Gamma, X})$ is an isomorphism, then restricting to $X - L $, we find ${ \mathcal O_X}  \lto \pi_{*}( {\mathcal O_S})  $ is an isomorphism there, so $S- \pi^{-1}(L) \lto X -L$ is an isomorphism. 
}
\end{Remark}

\bs

\begin{Proposition}\label{L2} If $\pi: S \rightarrow X$ is a finite morphism of schemes, then the natural map $$H^1(X,  \pi_*(\mathcal O_{S}^{*}) ) \longrightarrow H^1(S,  \mathcal O_S^*)$$ is an isomorphism. 
\end{Proposition}

\demo First we will show that the first higher direct image sheaf $R^1 \pi_*(\mathcal O_{S}^{*})$ is zero. This sheaf is the sheaf associated to the presheaf which to each open subset $V$ in $X$ associates the group $H^1(\pi^{-1}(V),\mathcal O^*_S|_{\pi^{-1}V})$ \cite[III, 8.1]{HAG}. Hence the stalk of this sheaf at a point $x\in X$ is the direct limit 
$$
\varinjlim_{x\in V} H^1(\pi^{-1}(V),\mathcal O^*_S|_{\pi^{-1}V}) \, .$$
An element in this direct limit is represented by a pair $(V, \alpha)$ where $V$ is an open set of $X$ containing $x$ and $\alpha \in  H^1(\pi^{-1}(V),\mathcal O^*_S|_{\pi^{-1}V})$. This group is just ${\rm Pic}(\pi^{-1}V)$, so the element $\alpha$ corresponds to an invertible sheaf $\mathcal L$ on $\pi^{-1}V$. We may assume that $V$ is affine, since affine open sets form a basis for the topology. Therefore, since $\pi$ is finite, the open subset $\pi^{-1}V$ of $S$ is also affine, and hence $\mathcal L$ is generated by global sections. Let $z_1, \ldots, z_r \in \pi^{-1}V$ be the finite set of points in $\pi^{-1}(x)$. We can find a section $s \in H^0(\pi^{-1}V, \mathcal L)$ that does not vanish at any of the  $z_1, \ldots, z_r$. So the zero set of  $s$ is a divisor $D$ whose support does not contain any of the $z_i$. Since $\pi$ is finite, it is a proper morphism, so $\pi(D)$ is closed in $V$ and does not contain $x$. Let $V'=V-\pi(D)$. Then ${\mathcal L}|_{\pi^{-1}(V')}$ is free, and since $\pi^{-1}(V') \subset \pi^{-1}(V)$, the image of $\alpha$ in the above direct limit is zero. Hence  $R^1 \pi_*(\mathcal O_{S}^{*})=0$. 
Now the statement of the lemma follows from the exact sequence of terms of low degree of the Leray spectral sequence
$$ 0\rightarrow H^1(X, \pi_*\mathcal O_S^*) \longrightarrow H^1(S, \mathcal O_S^*) \longrightarrow H^0(X, R^1 \pi_*(\mathcal O_{S}^{*})=0
$$
\QED

\ms

\begin{Note}{\rm Since $H^1(S, \mathcal O_S^*)$ computes ${\rm Pic} \, S$ see \cite[III, Ex. 4.5]{HAG}, this means that we can also compute ${\rm Pic}\, S$ as  $H^1(X, \pi_*\mathcal O_S^*)$. }
\end{Note}



\bs
\bs

\section{A local isomorphism for $\Apic$}

\ms 

In this section we prove a fundamental local isomorphism that allows us to compute the $\Apic$ group of a surface locally in terms of Cartier divisors on the curves $L$ and $\Gamma$. We first observe that if 
 $A$ is a local ring of dimension two satisfying $G_1$ and  $S_2$ with spectrum $X$ and punctured spectrum $X'$ then ${\rm A Pic} \, X={\rm Pic}\, X'$. Indeed ${\rm A Pic} \, X= {\rm A Pic} \,X'$ (see \cite[1.12]{H}), and $X'$ has no points of codimension two, so  ${\rm A Pic} \,X'= {\rm Pic}\, X'$.

\ms 

\begin{Theorem}\label{TL} Adopt assumptions \ref{A1}. Further assume that $X$ is the spectrum of a two dimensional local ring, $S$ is smooth,  and the map induced by $\pi$ from $\mathcal I_{L,X}$ to $  \pi_{*} (\mathcal I_{\Gamma, X})$ is an isomorphism. Then the map $$ \varphi: {\rm A Pic} X \longrightarrow {\rm Cart} \, \Gamma/ \pi^{*}{\rm Cart} \, L$$ in Proposition \ref{P2} is an isomorphism. 
\end{Theorem}
\demo  Let $x$ be the closed point of $X$. Set $X'=X - \{x\}$ and $S'=S - \{\pi^{-1}(x)\}$. 
As we noted above we can calculate ${\rm A Pic} \, X$ as ${\rm Pic} \, X'$ which is also $H^1(X', \mathcal O_{X'}^*)$. We consider sheaves of abelian groups on $X$
\begin{equation}\label{E2} 0 \rightarrow \mathcal O_{X}^{*}  \longrightarrow  \pi_*(\mathcal O_{S}^{*})   \longrightarrow
 \mathcal N  \rightarrow 0
\end{equation}
and similarly with primes
\begin{equation}\label{E3} 0 \rightarrow \mathcal O_{X'}^{*}  \longrightarrow  \pi_*(\mathcal O_{S'}^{*})   \longrightarrow
 \mathcal N'  \rightarrow 0
\end{equation}
Computing cohomology on $X$ along (\ref{E2}) we obtain the  exact sequence:
$$
0 \rightarrow   H^0(X, \mathcal O_{X}^{*})   \longrightarrow   H^0( X, \pi_*(\mathcal O_{S}^{*}) )  \longrightarrow
 H^0( X,\mathcal N)  \longrightarrow  {\rm Pic}\, X=0
$$
where $H^1(X, \mathcal O_{X}^*)={\rm Pic}\, X=0$ since $X$ is a local affine scheme. Now computing cohomology on $X'$ along
(\ref{E3})  we obtain the exact sequence
$$0 \rightarrow   H^0(X', \mathcal O_{X'}^{*})    \rightarrow   H^0( X', \pi_*(\mathcal O_{S'}^{*}) ) \rightarrow 
 H^0( X',\mathcal N')   \rightarrow {\rm A Pic} X \rightarrow     H^1(X',  \pi_*(\mathcal O_{S'}^{*}) )=0
$$
where $H^1(X',  \mathcal O_{X'}^{*} )={\rm A Pic} X$ and $H^1(X',  \pi_*(\mathcal O_{S'}^{*}) )=0$ because  by the Proposition \ref{L2} we have  $H^1(X',  \pi_*(\mathcal O_{S'}^{*}) )=H^1(S',  \mathcal O_{S'}^{*} )$. Now the latter is  ${\rm Pic}\, S'$, which in turn is equal to ${\rm A Pic} \, S$.  
But $ {\rm A Pic} \, S={\rm Pic} \, S$ because $S$ is smooth and finally ${\rm Pic} \, S=0$ because $S$ is a semi-local affine scheme.

Since $X$ and $S$ both satisfy $S_2$ any section of $\mathcal O_X$ or $\mathcal O_S$ over $X'$ or $S'$ extends to all of $X$ or $S$. Thus 
 $H^0(X, \mathcal O_{X}^{*}) = H^0(X', \mathcal O_{X'}^{*})$ and $H^0(X, \mathcal O_{S}^{*}) = H^0(X', \mathcal O_{S'}^{*})$. This allows us  to combine the above two sequences of cohomology into one:
 \begin{equation}\label{E4}0 \rightarrow H^0( X,\mathcal N)  \longrightarrow  H^0( X',\mathcal N')   \longrightarrow {\rm A Pic} X\ \longrightarrow 0
 \end{equation}
 By Proposition \ref{2.6} applied to both maps $S \lto X$ and $S' \lto X'$, we obtain $H^0( X,\mathcal N)=H^0( X,\mathcal N_0)$ and $H^0( X',\mathcal N')=H^0( X',\mathcal N'_0)$. Thus we  turn (\ref{E4}) into the following short exact sequence
   \begin{equation*} 0 \rightarrow H^0( X,\mathcal N_0)  \longrightarrow  H^0( X',\mathcal N'_0)   \longrightarrow {\rm A Pic} X\ \longrightarrow 0
 \end{equation*}
 Using this exact sequence we can derive the following diagram:
 \begin{equation*}
\begin{array}{rccccccc}
& 0  &  &  0 & & 0 & &  \\
& \downarrow &  & \downarrow &  & \downarrow &  &  \\
0 \rightarrow & H^0( X,\mathcal O_{L}^{*} )   & \longrightarrow & H^0( X,\pi_*(\mathcal O_{\Gamma}^{*})  )   & \longrightarrow
&  H^0( X,\mathcal N_0) & \rightarrow & 0 \\
 &\downarrow &  & \downarrow &  & \downarrow &  &  \\
0 \rightarrow & H^0( X',\mathcal O_{L'}^{*} )  & \longrightarrow & H^0( X',\pi_*(\mathcal O_{\Gamma'}^{*}) )  & \longrightarrow
&  H^0( X',\mathcal N'_0)  & \rightarrow & 0 \\
& \downarrow &  & \downarrow &  & \downarrow &  &  \\
 0 \rightarrow  & {\rm Cart} \, L & \longrightarrow & {\rm Cart}\,  \Gamma & \longrightarrow
&  {\rm A Pic} X & \rightarrow & 0 \\
 & \downarrow &  & \downarrow &  & \downarrow &  &  \\
& 0  &  &  0 & & 0 & &  \\
\end{array}
\end{equation*}
The first two rows in the  diagram are obtained applying cohomology to the short exact sequences 
$$0 \rightarrow \mathcal O_{L}^{*}  \longrightarrow  \pi_*(\mathcal O_{\Gamma}^{*})   \longrightarrow
 \mathcal N_0  \rightarrow 0$$
 and 
 $$0 \rightarrow \mathcal O_{L'}^{*}  \longrightarrow  \pi_*(\mathcal O_{\Gamma'}^{*})   \longrightarrow
 \mathcal N'_0  \rightarrow 0$$
and observing that again $H^1(X, \mathcal O_L^*)={\rm Pic}\,  L=0$ because $L$ is a local affine scheme, and $H^1(X', \mathcal O_{L'}^{*} )={\rm Pic} \, L'=0$ because $L'$ as a scheme is a disjoint union of generic points. 
The vertical columns arise from the fact that $L$ and $\Gamma$ are (semi)-local curves, so that when we remove the closed points we obtain the local rings of the generic points, namely the total quotient rings of $\mathcal O_L$ and $\mathcal O_{\Gamma}$, and the Cartier divisors are nothing else than the quotients of the units in the total quotient rings divided by the  units of the (semi)-local rings, i.e.    ${\rm Cart} \, L= K_{L}^*/\mathcal O_{L}^{*} $ and ${\rm Cart} \, \Gamma= K_{\Gamma}^*/\mathcal O_{\Gamma}^{*} $. Now the Snake Lemma yields the diagram. The last row of the above diagram implies the desired statement, namely $${\rm A Pic}\, X \cong {\rm Cart} \, \Gamma / \pi^*{\rm Cart}\, L.$$  \QED
\bs

\begin{Remark}\label{R3}{\rm  If $S$ is not smooth, we can replace ${\rm A Pic} \, X$ with the group $G$ defined in Remark \ref{2.5}, in which case the proof of Theorem \ref{TL} shows that $\varphi: G \longrightarrow {\rm Cart} \, \Gamma/ \pi^{*}{\rm Cart} \, L$ is an isomorphism, because the cokernel of $H^0({\mathcal N}_0) \longrightarrow H^0({\mathcal N}^{'}_0)$ is just the kernel of ${\rm A Pic}\, X \to {\rm A Pic}\, S$, and in our local case, for the image an element of ${\rm A Pic}\, X$ to vanish in ${\rm A Pic}\, S$ is the same as saying it is locally free, hence it is in ${\rm Pic}\, S$, which is zero because $S$ is a semi local ring.   }
\end{Remark}

\ms

\begin{Proposition}\label{completion}
If $A$ is a local ring of dimension one satisfying $S_1$ and $\cA$ is its completion, then $${\rm Cart}\, A= {\rm Cart} \, \cA$$
\end{Proposition}
\demo This follows for instance from the proof of \cite[2.14]{H} where it is shown that  ${\mathfrak a} \mapsto \cma$ gives a one-to-one correspondence between ideals of finite colength of $A$ and $\cA$, under which principal ideals corresponds to principal ideals.
\QED

\ms

The following proposition shows that the local calculation of $\Apic$ depends only on the analytic isomorphism class of the singularity when the normalization is smooth.



\begin{Proposition}\label{6.1} Let $A$ be a reduced two dimensional  local ring  satisfying $G_1$ and $S_2$ whose normalization $\nA$ is regular. 
Then 
$$
\Apic(\Spec A)=\Apic(\Spec \cA).
$$
\end{Proposition}
\demo
We let $S$ be the normalization of $X=\Spec A,$ take $L$ to be the conductor and $\Gamma$ to be $\pi^{-1}(L)$. Then by Theorem \ref{TL} (cf. Remark \ref{R2}) we can compute $\Apic X=\Cart \Gamma/\pi^* \Cart L$. On the other hand, we have shown (Proposition \ref{completion}) that the Cartier divisors of a one-dimensional local ring are the same as those of its completion. So applying Theorem \ref{TL} also to $\Spec \cA$ we prove the assertion.
\QED
\bs

The next example shows that we cannot drop the assumption on $\nA$ being regular.

\begin{Example}{\rm The divisor class groups of normal local rings and their behavior under completion have been studied by Mumford, Samuel, Scheja, Brieskorn, and others. In particular, if  $A = \CC[x,y,z]/(x^r + y^s + z^t)$ localized at the maximal ideal $(x,y,z)$, where $r < s < t$ are pairwise relatively prime, then $A$ is a unique factorization domain, so $\Apic(\Spec A)$ is 0. However, the completion $\nA$ of $A$ is not a unique factorization domain, so $\Apic (\Spec \nA)$ is not 0, except in the unique case $(r,s,t) = (2,3,5)$ \cite{S}.
}
\end{Example}

\bs

\section{Global exact sequences}

\ms

In this section we compare $\Pic$ and $\Apic$ of any surface $X$ to its normalization. This generalizes \cite[6.3]{H} which dealt with the case of a ruled cubic surface. In particular our result applies to a surface with ordinary singularities whose normalization is smooth, thus providing an answer to the hope expressed in \cite[6.3.1]{H}.

\begin{Theorem} \label{TG}Adopt assumptions \ref{A1}. Further assume that $X$ is a surface either affine or projective
and the map induced by $\pi$ from $\mathcal I_{L,X}$ to $  \pi_{*} (\mathcal I_{\Gamma, X})$ is an isomorphism. Then there is an  exact sequence:
\[  \hspace{-1.7in}(a)\qquad  {\rm Pic} \, X \longrightarrow {\rm Pic} \, S \longrightarrow {\rm Pic} \, \Gamma/ \pi^* {\rm Pic}\, L.
\]
\noindent
Furthermore, if $S$ is smooth, then  there is also an exact sequence $$ (b) \qquad {\rm A Pic} \, X \longrightarrow {\rm Pic} \, S \oplus  {\rm Cart} \, \Gamma/ \pi^* {\rm Cart}\, L\longrightarrow {\rm Pic} \, \Gamma/ \pi^* {\rm Pic}\, L \rightarrow 0.$$ 

\end{Theorem}
\demo For $(a)$ we use the natural map from Proposition \ref{P1}$(a)$ of ${\rm Pic}\, X$ to ${\rm Pic}\, S$ and the restriction map of ${\rm Pic}\, S$ to $ {\rm Pic} \, \Gamma/ \pi^* {\rm Pic}\, L$.  The composition is clearly zero, since a divisor class originating on $X$ will land in $ \pi^* {\rm Pic}\, L$. To show exactness in the middle, we recall the result of Proposition \ref{2.6} which shows that the sheaves $\mathcal N$, ${\mathcal N}_0$ in the following diagram are isomorphic:
\begin{equation*}
\begin{array}{ccccccccc}
0 &\rightarrow & \mathcal O_{X}^{*}  & \longrightarrow &   \pi_*\mathcal O_{S}^{*}  & \longrightarrow
& \mathcal N & \rightarrow & 0 \\
& & \scriptstyle{\alpha}\downarrow &  & \scriptstyle{\beta}\downarrow  &  &\scriptstyle{\gamma}\downarrow\scriptstyle{\cong}  &  &  \\
0 &\rightarrow & \mathcal O_{L}^{*} & \longrightarrow & \pi_* \mathcal O_{\Gamma}^{*} & \longrightarrow
& \mathcal N_0 & \rightarrow & 0 \\
\end{array}
\end{equation*}

Taking cohomology on $X$ and using Proposition \ref{L2}, we obtain a diagram of exact sequences
\begin{equation*}
\begin{array}{ccccccccccccc}
0 &\rightarrow & H^0(\mathcal O_{X}^{*})  & \rightarrow &  H^0( \pi_*\mathcal O_{S}^{*})  & \rightarrow &H^0(\mathcal N) &\rightarrow & {\rm Pic}\, X & \rightarrow  & {\rm Pic}\, S & \rightarrow  & H^1(\mathcal N)\\
& & \downarrow &  & \downarrow  &  &\downarrow\scriptstyle{\cong}  &  &  \downarrow &  & \downarrow  &  &\downarrow\scriptstyle{\cong}  \\
0 &\rightarrow & H^0(\mathcal O_{L}^{*})  & \rightarrow  &  H^0( \pi_*\mathcal O_{\Gamma}^{*})  &\rightarrow 
&H^0(\mathcal N_0) & \rightarrow  & {\rm Pic}\, L & \rightarrow  & {\rm Pic}\, \Gamma & \rightarrow  & H^1(\mathcal N_0)\\\end{array}
\end{equation*}
\noindent
From this sequence, we see that if an element of ${\rm Pic}\, S$ becomes zero in ${\rm Pic} \, \Gamma/ \pi^* {\rm Pic}\, L$, then it is zero in $H^1(\mathcal N_0)=H^1(\mathcal N)$ hence comes from an element of ${\rm Pic}\, X$.

For $(b)$ we first define the maps involved in the sequence. We use the natural map from Proposition \ref{P1}$(b)$ of ${\rm A Pic}\, X$ to ${\rm A Pic}\, S$, which is equal to ${\rm Pic}\, S$ since $S$ is smooth, together with the map $\varphi$ of Proposition \ref{P2}  applied locally. Note that $ {\rm Cart} \, \Gamma/ \pi^* {\rm Cart}\, L$ is simply the direct sum of all its contributions at each one of its points, since $L$ and $\Gamma$ are curves. 
The second map of $(b)$  is composed of the map ${\rm Pic} \, S \rightarrow {\rm Pic} \, \Gamma/ \pi^* {\rm Pic}\, L$ of $(a)$ and the natural maps of Cartier divisors to ${\rm Pic}$. 

The composition of the two maps is zero, because if we start with something in ${\rm A Pic}\, X$, then according to the construction of Proposition \ref{P2}, its two images in $ {\rm Pic} \, \Gamma/ \pi^* {\rm Pic}\, L$ will be the same. The second map of the sequence $(b)$ is clearly surjective.

To show exactness in the middle of $(b)$, suppose that a class $\d \in {\rm Pic}\, S$ and a divisor $D$ in ${\rm Cart}\, \Gamma$ have the same image in $ {\rm Pic} \, \Gamma/ \pi^* {\rm Pic}\, L$. First, we can modify $D$ by an element of  $\pi^* {\rm Cart}\, L$ so that $\d$ and $D$ will have the same image in $ {\rm Pic} \, \Gamma$. Next, by adding some effective divisors linearly equivalent to $mH$, where $H=0$ in the affine case and $H$ is a hyperplane section in the projective case, we can reduce to the case where $\d$ and $D$ are effective. 
Consider the exact sequence of sheaves 
$$ 0 \rightarrow {\mathcal I}_{\Gamma}(\d) \longrightarrow  {\mathcal O}_S(\d)\longrightarrow {\mathcal O}_{\Gamma}(\d) \rightarrow 0$$
If $X$ is affine, then $H^0(X, {\mathcal O}_S(\d))\rightarrow H^0(X, {\mathcal O}_{\Gamma}(\d))$ is surjective, so the section $s \in H^0(X, {\mathcal O}_{\Gamma}(\d))$ defining the divisor $D$ will lift to a section $s' \in H^0(X, {\mathcal O}_S(\d))$. This section defines a curve $\mathcal C$ in $S$ in the divisor class $\d$, not containing any component of $\Gamma$ in its support, whose intersection with $\Gamma$ is $D$. We can transport $\mathcal C$ restricted to $S-\Gamma$ to $X - L$, and take its closure in $X$. This will be an element of ${\rm A Pic}\, X$ giving rise to the $\d$ and $D$ we started with.

If $X$ is projective, we use a hyperplane section $H$. Since $H$ comes from ${\rm Pic}\, X$, it is sufficient to prove the result for $\d+mH$ and $D+mH$. Now for $m \gg 0$ the cokernel of the map $H^0(X, {\mathcal O}_S(\d+mH))\rightarrow H^0(X, {\mathcal O}_{\Gamma}(\d+mH))$ lands in $H^1(X,  {\mathcal I}_{\Gamma}(\d+mH))$, which is zero by Serre's vanishing theorem. Then the proof proceeds as in the affine case. 
\QED

\begin{Proposition}\label{injective}
With the hypotheses of Theorem \ref{TG} the following conditions are equivalent:
\begin{enumerate}
\item[(i)] the map  ${\rm Pic} \, X \longrightarrow {\rm Pic} \, S$  is injective
\item[(ii)] the map ${\rm Pic} \, L \rightarrow {\rm Pic}\, \Gamma$ is injective and $${\rm coker} (H^0(\mathcal O_{X}^{*} ) \to \pi_*H^0(\mathcal O_{S}^{*} ))={\rm coker} (H^0(\mathcal O_{L}^{*} ) \to \pi_*H^0(\mathcal O_{\Gamma}^{*} )).$$
\end{enumerate}
Furthermore if $S$ is smooth, conditions $(i)$ and $(ii)$ are also equivalent to 
\begin{enumerate}
\item[(iii)] The first map of Theorem \ref{TG}(b) is  injective
\end{enumerate}

\noindent
In addition, without assuming $S$ smooth, if conditions $(i)$ and $(ii)$ hold and the map  ${\rm Pic} \, L \rightarrow {\rm Pic}\, \Gamma$ is an isomorphism then  ${\rm Pic} \, X \rightarrow {\rm Pic}\, S$ is also an isomorphism.
\end{Proposition}
\demo
From the diagram of exact sequences in the proof of Theorem \ref{TG}(a), statement $(i)$ is equivalent to the exactness of the sequence:
\begin{equation}\label{E4.2}
0 \rightarrow  H^0(\mathcal O_{X}^{*})   \rightarrow H^0( \pi_*\mathcal O_{S}^{*})  \rightarrow H^0(\mathcal N) \to 0. 
\end{equation}
Since $H^0(X, \mathcal N)\cong H^0(\mathcal N_0) $, the exactness of  (\ref{E4.2}) implies the exactness of 
\begin{equation*}
0 \rightarrow  H^0(\mathcal O_{L}^{*})   \rightarrow H^0( \pi_*\mathcal O_{\Gamma}^{*})  \rightarrow H^0(\mathcal N_0) \to 0. 
\end{equation*}
Looking again at the diagram of exact sequences in the proof of Theorem \ref{TG}(a), this  implies $(ii)$. On the other hand, $(ii)$ clearly implies $(i)$. Now if $S$ is smooth, because of the local isomorphism of Theorem \ref{TL}, any element in the kernel of the first map of Theorem \ref{TG}(b) is zero in all the local groups ${\rm A Pic}\, ({\rm Spec}\, \mathcal O_{X,x})$, hence by Jaffe's sequence (see Proposition \ref{P1}(c)) is already in ${\rm Pic}\, X$. Thus $(iii)$ is also equivalent to $(i)$ and $(ii)$. 

The last statement follows again from the diagram of exact sequences in the proof of Theorem \ref{TG}(a).
\QED

\begin{Remark}{\rm If $X$ is integral  and projective in Proposition \ref{injective} then so is $S$, thus  $H^0( X,\mathcal O_{X}^{*} )=H^0( X,\pi_*\mathcal O_{S}^{*} )=k^*$, where $k$ is the ground field. Therefore ${\rm coker} (H^0(\mathcal O_{X}^{*} ) \to \pi_*H^0(\mathcal O_{S}^{*} ))=0$ and the equality of the cokernels in $(ii)$ holds if and only if $$
H^0(\mathcal O_{L}^{*} )=H^0(\mathcal \pi_*O_{\Gamma}^{*} ).$$  }
\end{Remark}

\s

\begin{Remark}\label{4.4}{\rm If $S$ is not smooth, then as in Remarks \ref{2.5} and \ref{R3} we can obtain the same results as in Theorem \ref{TG}, with ${\rm A Pic}\, X$ replaced by $G$.}
\end{Remark}

\s

The next theorem shows, that at least over the complex number  $\CC$, the map $\Pic X \to \Pic S$ is always injective.

\begin{Theorem}\label{4.5} Let $X$ be an integral surface in $\mathbb P^3$ over $k=\CC$. If $S$ is the normalization of $X$, then the natural map
$$\Pic X \lto \Pic S$$ is injective.
\end{Theorem}
\demo From the exponential sequence \cite[Appendix B, \S5]{HAG} we obtain an exact sequence of cohomology
$$ 0\to H^1(X_h, \Z) \lto H^1(X, \O_X) \lto \Pic X \lto H^2(X_h, \Z) \lto H^2(X, \O_X) \lto \ldots
$$
where $X_h$ is the associated complex analytic space of $X$. Now $H^1(X, \O_X)=0$ since $X$ is a complete intersection variety of dimension at least 2 \cite[III, Ex. 5.5]{HAG}. Furthermore $H^2(X_h, \Z)$ is a finitely generated abelian group, so we conclude that $\Pic X$ is also a finitely generated abelian group. 

Next, using Grothendieck's method of comparing $\Pic X$ to the Picard group of the formal completion of $\mathbb P^3$ along $X$ \cite[IV, \S3]{HLN}, the proof of \cite[3.1]{HLN} and the groundfield $\CC$ being of characteristic zero shows that $\Pic X$ is torsion free \cite[Ex. 20.7]{HDT}. Thus $\Pic X$ is in fact a free finitely generated abelian group.

Taking $L$ to be the conductor and $\Gamma$ its inverse image in $S$, the Assumptions \ref{A1} are satisfied. Now, looking at the exact sequences used in the proof of Theorem \ref{TG}, since $X$ is integral and projective, $H^0(\O_X^*) \to H^0(\pi_*\O_S^*)$ is an isomorphism, so $H^0(\N)$ 
is equal to the kernel of the map $\Pic X \to \Pic S$. If it is non zero, it must be a finitely generated free abelian group. 

Since $L$ and $\Gamma$ are projective curves, 
the group of units in each is a direct sum of $k$-vector spaces and copies of $k^*$. To see this, refer to Proposition \ref{5.11}, and note that the first sequence splits, since $k^*$ is contained in $\O_{\C}^*$. Now comparing these sequences for $L$ and $\Gamma$ we see that the cokernel of the map $H^0(\O_L^*) \to H^0(\pi_* \O_{\Gamma}^*)$ must also be a direct sum of a $k$-vector space and copies of $k^*$. But $H^0(\N_0)=H^0(\N)$ is a finitely generated free abelian group, so this cokernel must be zero.

Hence $H^0(\N_0)$ is equal to the kernel of the map  $\Pic L \to \Pic \Gamma$. Since $L$ and $\Gamma$ are curves, there are degree maps on each irreducible component to copies of $\Z$. The map $\Gamma \to L$ is surjective, so an element of positive degree on $L$ remains an element of positive degree on $\Gamma$. Thus the kernel of the map $\Pic L \to \Pic \Gamma$ is just the kernel of the degree $0$ part $\Pic^0L \to \Pic^0\Gamma$. These are group schemes, successive extensions of abelian varieties by copies of $\G_a$ and $\G_m$ \cite{O1}. In particular, the kernel is also a group scheme, of finite type over $k$. Since $H^0(\N_0)$ is a finitely generated free abelian group, as a group scheme it must have dimension zero, and hence (again using characteristic zero) must be a finite abelian group. But it is also a free abelian group, hence it is zero. Thus $\Pic X \to \Pic S$ is injective. 
\QED



\section{Results on curves}
For our applications to surfaces, we need to know something about the curves $L$ and $\Gamma$. A {\it curve} in this section will be a one dimensional scheme without embedded points, hence satisfying the condition $S_1$ of Serre. For any ring $A$, we denote by ${\rm Cart}\, A$ the group of Cartier divisors of ${\rm Spec}\, A$.  We first compute the local group $\Cart A$ at a singular point of a curve in terms of the number of local branches and an invariant $\delta$. Then we study the Picard group of a projective curve showing the contribution of the singular points. The results of this section are essentially well-known (see the papers of Oort \cite{O} and \cite{O1} on the construction of the generalized Jacobian). Here we gather the results on groups of divisors and divisor classes that we will need later. 


\s

In the following $k^*$ is the multiplicative group of units of the field $k$ and $k^+$ is the whole field $k$ as a group under addition.
\begin{Theorem}\label{cart-curves}
Let $A$ be a reduced local ring of dimension one with residue field $k$ algebraically closed of characteristic zero.
Let $\nA$ be the normalization of $A$, let $\r$ be the number of maximal ideals of $\nA$ (the number of {\it branches} of the curve singularity), and let $\delta$ be the length of $\nA/A$. Then 
$$ {\rm Cart} \, A \cong \Z^{\r} \oplus (k^*)^{\r-1} \oplus (k^+)^{\delta-\r+1}.
$$
\end{Theorem}

\s

Before the proof we need some Lemmas.

\s

\begin{Lemma}\label{units} Let $A$ be a ring, $\a$ an ideal, and assume either $A$ is a local ring or $A$ is complete in the $\a$-adic topology. Then the natural map of units $A^* \to (A/\a)^*$ is surjective.
\end{Lemma}
\demo Let $a\in A$ be an element such that $\oa \in A/\a$ is a unit. In the local case this means $\oa \not\in \m_{A/\a}$, which is equivalent to saying $a\not\in \m_A$, so $a$ is a unit.

In the complete case, there exists a $b\in A$ such that $\oa\ob=\o1$. In other words, $ab=1+x$ for some $x\in \a$. Now let $u=1-x+x^2-x^3+\ldots$ which exists in $A$ since it is complete with respect to the $\a$-dic topology. Then $abu=1$, so $a$ is a unit. 
\QED

\bs

Recall that by $\a^+$ we denote the ideal $\a$ as a group under addition. 
\begin{Lemma}\label{exp} Let $A$ be a ring and $\a$ an ideal. Assume  $A$ is complete in the $\a$-adic topology and that $A$ contains the rational numbers $\Q$. 
Then there is an exact sequence of abelian groups 
$$ 0 \to \a^+ \stackrel{\alpha}\longrightarrow A^* \stackrel{\beta}\longrightarrow (A/\a)^* \to 1
$$
where the map $\alpha$ send an element $a$ of $\a$ to 
$${\rm exp}(a)= 1+a+ \frac{1}{2}a^2+\frac{1}{6}a^3+ \ldots$$ 
\end{Lemma}
\demo  
Clearly the composition $\beta\alpha$ is $1$, and we know that $\beta$ is surjective from the previous lemma. If $u$ is an element in the kernel of $\beta$ that means $u=1+y$ for some $y\in \a$, and then $$a={\rm log} \, u=y-\frac{1}{2}y^2+\frac{1}{3}y^3-\ldots$$
is in $\a$ and gives $u$ via the map $\alpha$. This shows also that $\alpha$ is injective, so the sequence is exact. We need only to know that the $\rm exp$ map and the $\rm log$ map are inverses to each other and this is purely formal.
\QED

\ms

{\bf{Proof of Theorem \ref{cart-curves}.} }
By Proposition \ref{completion} we may assume that $A$ is a complete local ring. Notice that the normalization of the completion is the completion of the normalization so $\nA$ is also complete with respect to its  Jacobson radical $J$. Let $K$ be the total quotient ring of $A$. Then $A \subset \nA \subset K$,  ${\rm Cart}\, A= K^*/A^*$ and ${\rm Cart}\, \nA= K^*/{\nA}^*$. Therefore, we obtain the exact sequence 
$$ 
0 \to 
{\nA}^*/ A^*
 \longrightarrow  {\rm Cart } \, 
A \longrightarrow {\rm Cart } \, \nA \to 0.
$$
Now $\nA$ is just a product of $\r$ discrete valuation rings, so ${\rm Cart}\, \nA = \Z^{\r}$ (recall that the group of Cartier divisors of a discrete valuation ring is $\Z$). Since $\Z$ is a free abelian group, the sequence splits and we can write 
$$ {\rm Cart} \, A \cong \Z^{\r} \oplus{\nA}^*/ A^*.
$$
To analyze this latter group we will apply Lemma \ref{exp} to the rings $A$ and $\nA$. Let $\m$ be the maximal ideal of $A$. Then we can write
\begin{equation}\label{D}
\begin{array}{rccccccc}
& 0  &  &  1 & & 1 & &  \\
& \downarrow &  & \downarrow &  & \downarrow &  &  \\
0 \rightarrow & \m   & \longrightarrow & A^*  & \longrightarrow
&  k^*& \rightarrow & 1 \\
 &\downarrow &  & \downarrow &  & \downarrow &  &  \\
0 \rightarrow & J  & \longrightarrow & \nA^*  & \longrightarrow
&  (k^*)^{\r} & \rightarrow & 1 \\
& \downarrow &  & \downarrow &  & \downarrow &  &  \\
 0 \rightarrow  & M & \longrightarrow &{\nA}^*/ A^*& \longrightarrow
&   (k^*)^{\r -1}  & \rightarrow & 1 \\
 & \downarrow &  & \downarrow &  & \downarrow &  &  \\
& 0  &  &  1 & & 1 & &  \\
\end{array}
\end{equation}
The first column of (\ref{D}) defines an $A$-module $M$ of finite length $\delta -\r+1$, as is evident from considering the analogue  diagram of ideals and rings without $^*$. The first two rows  of (\ref{D}) are the applications of Lemma \ref{exp} to $(A, \m)$ and to $(\nA, J)$. The bottom row of (\ref{D}) gives us an exact sequence for ${\nA}^*/ A^*$. This sequence splits because $A$, being complete of characteristic zero,  contains its residue field $k$. And finally, as an abelian group $M$ is just isomorphic to $(k^+)^{\delta-\r+1}$. This gives the desired decomposition.
\QED

\bs

In the next proposition we address the non-reduced case.

\begin{Proposition}\label{NRed} Let $A$ be a local ring of dimension one satisfying $S_1$ with residue field $k$ of characteristic zero, let $K$ be the total quotient ring of $A$, and let $I$ be the ideal of nilpotents of $A$.
Then there is an exact sequence
$$0 \to ( I \otimes K)/I \lto {\rm Cart} \, A \lto {\rm Cart}\, A_{\rm red} \to 0.
$$
\end{Proposition}

\demo
We apply Lemma \ref{exp} to the pairs $(A,I)$ and $(K, I \otimes K)$, obtaining a diagram
\begin{equation*}
\begin{array}{rccccccc}
& 0  &  &  1 & & 1 & &  \\
& \downarrow &  & \downarrow &  & \downarrow &  &  \\
0 \rightarrow & I   & \longrightarrow & A^*  & \longrightarrow
&   A_{\rm red}^*& \rightarrow & 1 \\
 &\downarrow &  & \downarrow &  & \downarrow &  &  \\
0 \rightarrow & I \otimes K  & \longrightarrow & K^*  & \longrightarrow
&   K_{\rm red}^*& \rightarrow & 1 \\
& \downarrow &  & \downarrow &  & \downarrow &  &  \\
 0 \to & (I \otimes K)/I & \lto& {\rm Cart} \, A & \lto & {\rm Cart}\, A_{\rm red} & \to 0\\
 & \downarrow &  & \downarrow &  & \downarrow &  &  \\
& 0  &  &  0 & & 0 & &  \\
\end{array}
\end{equation*}

\QED

\begin{Note} {\rm As an abelian group $(I \otimes K)/I $ is a $k$-vector space, in particular, when we assume that ${\rm char}\, k=0$, it is a torsion-free abelian group.}
\end{Note}

\bs

\begin{Examples}\label{5.5}{\rm
\begin{enumerate}
\item[(a)] If $A$ is a node, then $\r=2$ and $\delta=1$, so ${\rm Cart} \, A \cong \Z^{2} \oplus k^*$. This result can also be deduced computationally from \cite[3.1]{H}, since any principal ideal generated by a non-zerodivisor  in $A=k[[x,y]]/(xy)$ is of the form $(x^r+ay^s)$ with $a \in k^*$. Two of them multiply by adding the exponents of $x$ and $y$ placewise, and multiplying the coefficients of $a$.
\item[(b)] If $A$ is a cusp, then $\r=1$ and $\delta=1$, so ${\rm Cart} \, A \cong \Z \oplus k^+$. One could also recover this result from  \cite[3.7]{H}. 
\item[(c)] More generally, if $A$ is the local ring of a plane curve singularity, then $\r$ is the number of branches, and $\delta$ is the sum $\sum \frac{1}{2} r_i(r_i-1)$ taken over the multiplicities $r_i$ of the point itself and all the infinitely near singular points \cite[IV, Ex. 1.8 and V, 3.9.3]{HAG}. So for an ordinary plane triple point for example, $\r=3$, $\delta=3$, we have $ {\rm Cart} \, A \cong \Z^{3} \oplus (k^*)^{2} \oplus (k^+)$.
\item[(d)] Suppose $A$ is the local ring of a nonplanar triple point. Then we can suppose $$A=k[[x,y,z]]/(xy,xz,yz).$$ Generalizing the method of \cite[3.1]{H}, a principal ideal is generated by an element of the form $x^r+ay^s+bz^t$. Under multiplication the exponents of $x,y,z$ add, while the coefficients of $y,z$ multiply, respectively. Hence ${\rm Cart} \, A \cong \Z^{3} \oplus (k^*)^{2}$. From Theorem \ref{cart-curves} we infer that $\delta=2$. 
\item[(e)] Suppose that $A$ is the local ring of a point on a double line in the plane, such as $k[x,y]/(x^2)$ localized at $(x,y)$. Then $A_{\rm red}$ is $k[y]$ localized at $(y)$. The ideal of nilpotents $I$ is $xA$, which is a free module of rank $1$ over $A_{\rm red}$. Thus $(I \otimes K)/I\cong k(y)/k[y]\cong y^{-1}k[y^{-1}]$ is  an infinite dimensional $k$-vector space.  In the sequel we will denote this vector space by $W$.
\end{enumerate}
}
\end{Examples}

\bs 

Now let us study $H^0(X, \mathcal O_{\mathcal C}^*)$ and ${\rm Pic}\, \mathcal C$ for a projective curve $\mathcal C$. 

\begin{Proposition}\label{5.8} Suppose that $\mathcal C$ is an integral projective curve with normalization $\nC$. Then there is an exact sequence
$$ 0\to \pi_*\mathcal O_{\nC}^*/\mathcal O_{\C}^* \lto {\rm Pic}\, \C \lto {\rm Pic}\, \nC \to 0
$$

\end{Proposition}
\demo  
From the short exact sequence of sheaves of units 
$$1 \to \mathcal O_{\C}^* \lto  \pi_*\mathcal O_{\nC}^* \lto \pi_*\mathcal O_{\nC}^*/\mathcal O_{\C}^* \to 1
$$
taking cohomology, we obtain the long exact sequence
\begin{equation}\label{MN} 1 \to H^0(\mathcal O_{\C}^*) \lto  H^0(\pi_*\mathcal O_{\nC}^*)  \lto \pi_*\mathcal O_{\nC}^*/\mathcal O_{\C}^* \lto {\rm Pic}\, \C \lto {\rm Pic}\, \nC \lto 0\end{equation}
Indeed,  $ \pi_*\mathcal O_{\nC}^*/\mathcal O_{\C}^*$ is supported only at points, so $H^0(\pi_*\mathcal O_{\nC}^*/\mathcal O_{\C}^*)=\pi_*\mathcal O_{\nC}^*/\mathcal O_{\C}^*$ and it has no $H^1$. By Proposition \ref{L2}, $$H^1(\C,  \pi_*\mathcal O_{\nC}^*)=H^1(\nC,  \mathcal O_{\nC}^*)={\rm Pic}\, \nC.$$ Finally, since $\C$ and $\nC$ are both integral,   
$$H^0( \mathcal O_{\C}^*)= H^0(\pi_*\mathcal O_{\nC}^*)=k^*.$$ 
Thus we obtain the desired short exact sequence.
\QED 

\bs

\begin{Remark}\label{5.9}
{\rm  Since $\nC$ is a nonsingular projective curve, ${\rm Pic}\, \nC$  is an extension of $\Z$ by an abelian variety of dimension $g={\rm genus}\, \nC$. The dimension of ${\rm Pic}\, \C$ is given by $p_a(\C)$, the arithmetic genus of $\C$. On the other hand, the group $ \pi_*\mathcal O_{\nC}^*/\mathcal O_{\C}^*$ is supported at a finite number of points and can be computed as in Theorem \ref{cart-curves}, giving 
$$  \pi_*\mathcal O_{\nC}^*/\mathcal O_{\C}^*\cong (k^*)^{\sum (\r_i-1)}\oplus (k^+)^{\sum(\delta_i-\r_i+1)}
$$
where $\r_i$ and $\delta_i$ are defined at all the singular points of $\C$. Thus reading dimensions on the short exact sequence of Proposition \ref{5.8} we recover the well-known formula 
$$p_a(\C)=g(\nC) + \sum \delta_i.$$ 
The difference of the arithmetic genus and the geometric genus is $\delta$, which is  the total number of copies of $k^*$ or $k^+$. \cite[V, 3.9.2]{HAG} .  }\end{Remark}

\s



\begin{Proposition}\label{5.11}
 Suppose that $\C$ is non reduced projective curve. Let $\mathcal I$ be the sheaf of nilpotent elements of $\C$ and let $\C_{\rm red}$ be the reduced curve. Assume that ${\rm char} \, k =0$. Then there are two exact sequences

\begin{equation*}
\begin{array}{ccccccc}
0 \to & H^0(\mathcal I) &\lto & H^0(\mathcal O^*_{\C}) &\lto&  H^0(\mathcal O^*_{\C_{\rm red}})&\to 0 \\
\\
0 \to & H^1(\mathcal I) &\lto & {\rm Pic}\, \C &\lto&   {\rm Pic}\, \C_{\rm red} &\to 0 \\

\end{array}
\end{equation*}
\end{Proposition}
\demo
Just take cohomology of the sequence 
$$ 0 \to \mathcal I \lto \mathcal O^*_{\C} \lto \mathcal O^*_{\C_{\rm red}} \to 0
$$
where the first map is the exponential map defined as in Lemma \ref{exp}.  Since $\mathcal I$ is a coherent sheaf on a scheme of dimension one, $H^2(\mathcal I)=0$.
On each connected  component of $\C_{\rm red}$ the sections of $ \mathcal O^*_{\C_{\rm red}}$ are just $k^*$. These sections lift to $H^0(\mathcal O^*_{\C})$ so the long exact sequence splits in two, and the first sequence also splits.
\QED
\s


\begin{Example} \label{5.13}{\rm 
\begin{enumerate}
\item[(a)] If $\C$ is two lines in the plane $\mathbb P^2$ meeting at a point, then $\nC$ is two disjoint lines. We find that $H^0(\O_\C^*)=k^*$ and $H^0(\O_{\nC}^*)=(k^*)^2$. By Example \ref{5.5}(a) the contribution $\pi_*\O_{\nC}^*/\O_\C^*$ of the point is $k^*$, and clearly ${\rm Pic}\, \nC=\Z^2$. The exact sequence (\ref{MN}) then becomes
$$ 1 \to k^* \lto (k^*)^2 \lto k^* \lto \Pic \C \lto \Pic \nC \to 0
$$
hence ${\rm Pic} \, \C= \Z \oplus \Z$ (see also \cite[3.6]{H}).

\item[(b)] If $\C$ is a plane nodal cubic curve, then $\nC$ is a a rational curve with $\Pic \nC=\Z $. The curves are integral so we can apply Proposition \ref{5.8}. The local contribution of the node (see Example \ref{5.5}(a)) is $k^*$, so $\Pic \C=\Z\oplus k^*$ (see also \cite[II, Ex. 6.7]{HAG}).

\item[©] If $\C$ is a plane cuspidal cubic curve, the local contribution of the cusp is $k^+$ (see Example \ref{5.5}(b)), so as in (b) we find $\Pic \C=\Z\oplus k^+$ (see also \cite[II, 6.11.4]{HAG}). 

\item[(d)] If $\C$ is the union of three lines in the plane meeting at a single point $P$, then the local contribution is $(k^*)^2\oplus k^+$ (see Example \ref{5.5}(c)), and $\nC$ is three lines, so the exact sequence (\ref{MN}) is
$$
1 \to k^* \lto (k^*)^3 \lto (k^*)^2 \oplus k^+ \lto \Pic C \lto \Pic \nC \to 0
$$
Therefore $\Pic \C=\Z^3 \oplus k^+$.

\item[(e)] If $\C$ is the union of three lines in the plane making a triangle, then we have three local contribution of $k^*$ from the thee nodes, and using the sequence analogous to the one in $(d)$ above, we find $\Pic \C=\Z^3 \oplus k^*$. 

\item[(f)] Now suppose that $\C$ is the union of three lines in $\mathbb P^3$ meeting at a point but not lying on a plane. In this case the local contribution of the point is just $(k^*)^2$ as we can see from Example \ref{5.5}(d). Thus as in $(d)$ and $(e)$ we compute $\Pic \C=\Z^3$.

\end{enumerate}
}
\end{Example}

\bs

\begin{Remark}\label{5.14}{\rm 
Note that in all the examples above, the `dimension' of $\Pic \C$, meaning the number of factors $k^*$ or $k^+$ plus the dimension of the abelian variety (not present in these examples) is equal to the arithmetic genus $p_a$  of the curve. The plane cubic curves have $p_a$ equal one, while the non-planar cubic curve of $(f)$ has  $p_a$ equal zero.

}\end{Remark}


\section{Examples and Applications} 

\bs


Throughout this section, $X$ will denote a reduced surface in $\mathbb P^3$ and $\pi: S \to X$ will be its normalization. We take $L$ to be the singular locus of $X$ and $\Gamma=\pi^{-1}(L)$. Then the hypotheses of Assumptions \ref{A1} are satisfied. In fact $X$, being a hypersurface in $\mathbb P^3$, is Gorenstein, and so satisfies $G_1$ and $S_2$. 
We will compute $L$, $\Gamma$, $\Pic X$, $\Pic S$, and $\Apic ({\Spec} \O_{X,P})$ for some selected surfaces $X$ in $\mathbb P^3$, to illustrate the previous theoretical material.
\bs

\begin{Example}\label{two planes}{\rm 
 Let $X=H_1 \cup H_2$ be the union of two planes in $\mathbb P^3$ meeting along a line $L$. Then $S$ is the disjoint union of the two planes $H_1 \sqcup H_2$, and $\Gamma$ is a line in each of the $H_i$. The reduced line $L$ is the conductor of the normalization. So for a point $x \in L$ there will be two points of $\Gamma$ lying over it, and using Theorem \ref{TL} we can compute $$\Apic(\Spec \O_{X,x}) =\Cart \Gamma/\pi^* \Cart L= (\Z \oplus \Z)/\Z \cong \Z.
$$
To calculate $\Pic X$, since $X$ and $L$ have each one connected component, and $S$ and $\Gamma$ have each two connected components, we have ${\rm coker} (H^0(\mathcal O_{X}^{*} ) \to \pi_*H^0(\mathcal O_{S}^{*} ))={\rm coker} (H^0(\mathcal O_{L}^{*} ) \to \pi_*H^0(\mathcal O_{\Gamma}^{*} ))$. Furthermore $\Pic L= \Z \to \Pic \Gamma=\Z^2$ is injective, so by Proposition \ref{injective}, $\Pic X \to \Pic S$ is injective. 


 
Now by Theorem \ref{TG}, since $\Pic \Gamma/\pi^* \Pic L \cong \Z$ and $\Pic S =\Z \oplus \Z$, we conclude that $\Pic X=\Z$. Since  $\Apic(\Spec \O_{X,x})=\Z$ for each point $x\in L$, the sequence of  Proposition \ref{P1} $($c$)$ becomes \begin{equation}\label{J}
0 \to  \Pic X \lto \Apic X \lto \Div L \to 0,
\end{equation}
where $\Div L$ is just the direct sum of a copy of $\Z$ at each point of $X$. Note the last map is surjective, because using a sum of lines in one of the planes $H_i$ we can get any divisor on $L$. Thus the sequence (\ref{J}) splits and we obtain $\Apic X \cong \Z \oplus \Div L$. ( see \cite[5.2, 5.3, 5.4]{H} and the discussion following for more details about curves on $X$.)
}
\end{Example}


\s

\begin{Example}\label{three planes}{\rm Let $X$ be the union of three planes in $\mathbb P^3$ meeting along a line $L_0$. Then $S$ is the disjoint union of the three planes and $\Gamma_0$ is the union of one line in each of the planes. Taking $L_0$ and $\Gamma_0$ to be reduced, for a point $P \in L_0$ we have $\Cart_P L_0=\Z$ and $\Cart_{Q_i} \Gamma_0=\Z$ for each of the three points $Q_1, Q_2, Q_3$ lying over $P$. So the quotient $\Cart_Q \Gamma_0/\pi^* \Cart_P L_0 $ that appears in Proposition  \ref{P2} is just $\Z^2$. 

However in this case the conductor of the normalization is not $L_0$ reduced! It is the scheme structure supported on the line $L_0$ defined by $\I_{L_0, {\mathbb P^3}}^2$ ( see Example \ref{5.5}$($c$)$, which shows that $\delta=3$ for an ordinary plane triple point).

Lifting $L$ up to $S$, then $\Gamma$ consists of a planar double line in each plane. Now for $P\in L$, the ideal of nilpotents $I$ in the local ring $k[x,y,z]_{(x,y,z)}/(x,z)^2$ is free of rank $2$ over $k[y]$, so using Example \ref{5.5}(e) we obtain $\Cart_P L =\Z \oplus W^2$ (where $W$ is the $k$-vector space $y^{-1}k[y^{-1}]$) and at each point $Q_i \in \Gamma$ above $P$, we find $\Cart_{Q_i} \Gamma= \Z \oplus W.$ The quotient $\Cart_Q \Gamma/\pi^* \Cart_P L$ is $\Z^2 \oplus W$. This, therefore is the group $\Apic (\Spec \O_{X,P})$, not just the $\Z^2$ found above.

Note that if $\I$ is the sheaf of nilpotents of $\O_L$ (the non-reduced structure), then $\I \cong \O_{L_0}(-1)^2$. This sheaf has no cohomology on the line $L_0$, so $H^0(\O^*_L)=k^*$ and $\Pic L = \Z$. Similarly, for each component $\Gamma_{i}$ of $\Gamma$ we have  $H^0(\O_{\Gamma_{i}}^*)=k^*$ and $\Pic \Gamma_{i}=\Z$. Now from Proposition \ref{injective} we conclude that $\Pic X \lto \Pic S$ is injective. Looking at the first sequence of Theorem \ref{TG} we obtain that $\Pic X =\Z$.
}
\end{Example}

\s

\begin{Example}\label{6.6}$[$Pinch points and the ruled cubic surface$]$ 
{\rm A {\em pinch point} of a surface is a singular point that is analytically isomorphic to $k[x,y,z]_{(x,y,z)}/(x^2z-y^2)$. A typical example is the ruled cubic surface
$X$ in \cite[\S 6]{HAG}. This surface has a double line $L$ with two pinch points. The normalization $S$ is a nonsingular ruled cubic surface (scroll) in $\mathbb P^3$. The inverse of $L$ is a conic $\Gamma$ in $S$, and the restriction $\pi: \Gamma \to L$ is a 2-1 mapping, ramified over the two pinch points. Thus at a pinch point  $x \in L$, there is just one point $z\in \Gamma$ above $x$, and the mapping $\Cart L =\Z \to \Cart \Gamma=\Z$ is multiplication by $2$. Thus $\Apic (\Spec \O_{X,x}) =\Z/2\Z$ at the pinch point. 

Since $X$ is integral, $L$ and $\Gamma$ are both integral, and $\Pic L \to \Pic \Gamma$ is injective,  Proposition \ref{injective} applies and we find $\Pic X \to \Pic S$ is injective. As $\Pic \Gamma/ \pi^* \Pic L= \Z/2\Z$, and $\Pic S=\Z \oplus \Z$, we obtain also $\Pic X \cong \Z \oplus \Z$. The global sequences from Theorem \ref{TG} now give the results of  \cite[6.1, 6.2, 6.3]{H} with simplified proofs. 
}
\end{Example}




\begin{Example}\label{Steiner}$[$The Steiner surface, see \cite[page 137]{SR}, \cite[page 478]{P}, \cite[Ex. 5.5]{R}$]$ 
{\rm According to Pascal \cite{P}, this `superficie romana di Steiner' was discovered by Jakob Steiner in 1838 during a visit to Rome. It was first published by Kummer in 1863 in his article on quartic surfaces containing infinitely many conics, where he attributed it to Steiner. The Steiner surface $X$ is a projection of the Veronese surface $S$ in $\mathbb P^5$ (the 2-uple embedding of $\mathbb P^2$ in $\mathbb P^5$), and has the equation $$
x^2y^2+x^2z^2+y^2z^2=xyzw$$
in $\mathbb P^3$. To see this, let $t, u, v$ be coordinate in $\mathbb  P^2$, so that $S$ is given by $t^2, tu, u^2, tv, uv, v^2$ in $\mathbb P^5$, and project by taking
\begin{eqnarray*}
x&=&tu \\
y&=&tv\\
z&=&uv\\
w&=&t^2+u^2+v^2.
\end{eqnarray*}
Then $x, y, z, w$ satisfy the equation above. The singular locus $L$ of $X$ consists of the three lines $x=y=0$, $x=z=0$, $y=z=0$ meeting at the point $P=(0, 0, 0, 1)$ in $\mathbb P^3$. The curve $L$ is a double curve for $X$, with two pinch points on each line. It is the conductor of the normalization $\pi: S \to X$. The inverse image of $L$ in $S$ is a curve $\Gamma$, consisting of three conics, each meeting the other two in a point, which are the images of the three lines $t=0$, $u=0$, $v=0$ of $\mathbb P^2$, forming a triangle. The three vertices of the triangle go to the triple point $P$, while each side of the triangle is a 2-1 covering of the corresponding line.  

At the triple point $P$ we have $\Cart_P L= \Z^3 \oplus (k^*)^2$ by Example \ref{5.5}(d). On $\Gamma$, there are three nodes $Q_1$, $Q_2$, $Q_3$ lying over $P$, at each of which $\Cart_{Q_i} \Gamma =\Z^2 \oplus k^*$ by Example \ref{5.5}(a). Thus $\Apic (\Spec \O_{X,P})=(\Z^2 \oplus k^*)^3/\Z^3\oplus (k^*)^2=\Z^3 \oplus k^*$. At any other point of $L$ we have either an ordinary double point with $\Apic$ equal to $\Z$ (by Example \ref {two planes}) or a pinch point  with $\Apic$ equal to $\Z/2\Z$  (see Proposition \ref{6.1} and  Example \ref{6.6}). 

From Examples \ref{5.13}(f) and \ref{5.13}(e) we know that $\Pic L=\Z^3$ and $\Pic \Gamma= \Z^3 \oplus k^*$. The map $\Pic L \to \Pic \Gamma$ sends a generator of one line to twice the generator of the corresponding line in $\Gamma$. So $\Pic L \to \Pic \Gamma$ is injective and the quotient is $(\Z/2\Z)^3 \oplus k^*$.

We look at the first sequence of Theorem \ref{TG} 
$$
\Pic X \lto \Pic S \lto \Pic \Gamma/\pi^*\Pic L.$$
Since $S$ is the Veronese surface in $\mathbb P^5$, which is the 2-uple embedding of $\mathbb P^2$, we have $\Pic S=\Z$, and the hyperplane section $H$ is twice a generator. A curve $\C$ in $\Pic S$ goes to zero  in $\Pic \Gamma/\pi^*\Pic L$ if and only if its intersection with each line of $\Gamma$ is even, which means that $\C$ has even degree in $\P^2$, so in $\Pic S$ it is a multiple of $H$. Hence the kernel of $\Pic S \to \Pic \Gamma/\pi^* \Pic L$ is just $\Z\cdot H$. On the other hand, Proposition \ref{injective} applies to show that $\Pic X \to \Pic S$ is injective, thus $\Pic X=\Z$, generated by $H$.  (Note that since $X$ in $\mathbb P^3$ is a projection of $S$ in $\mathbb P^5$, the hyperplene class $H$ on $X$ lifts to the hyperplane class on $S$).

Finally, we can describe $\Apic X$  using Theorem \ref{TG}(b).
}
\end{Example}

\s

\begin{Definition}\label{new6.6} {\rm We say a reduced surface $X$ in $\P^3$ has {\it ordinary singularities} if its singular locus $L$ consists of a double line with transversal tangents  at most points, plus a finite number of pinch points and non-planar triple points.}
\end{Definition}
\s
\begin{Remark}\label{new6.7} {\rm The significance of ordinary singularities is that one knows from the literature, at least in characteristic zero, that the generic projection $X$ in $\P^3$ of any nonsingular surface $S$ in $\P^n$ has only ordinary singularities \cite{MP}. In characteristic $p>0$ the same applies after replacing $S$ if necessary by a suitable $d$-uple embedding \cite{R0}. Conversely, if $X$ in $\P^3$ has ordinary singularities its normalization $S$ is smooth (but may not have an embedding in $\P^n$ of which $X$ is the projection).}
\end{Remark}
\s
\begin{Proposition}\label{new6.8} If $X$ is a surface in $\P^3$ with ordinary singularities, then we can describe ${\rm APic} \, X$ by the sequence of Proposition \ref{P1}(3) as
$$
0 \to \Pic X \lto \Apic X \lto  \bigoplus_{P \in X} \Apic ({\rm Spec} \, \O_{X,P})$$
where
$$
\Apic ({\rm Spec} \, \O_{X,P})=\left\{ \begin{array}{cc}
\Z &\mbox{ {\rm at a general point of $L$}} \\
  \Z/2\Z &\mbox{ {\rm at a pinch point of $L$}}\\
  \Z^3 \oplus k^* &\mbox{{\rm at a triple point of $L$\, .}}
       \end{array} \right.
$$
\end{Proposition}
\demo Indeed, the local calculation of $\Apic$ is stable under passing to the completion (see Proposition \ref{6.1}) and the calculations for these three kinds of points have been done Examples \ref{two planes}, \ref{three planes} and \ref{6.6}.
\QED

\s
\begin{Example}\label{new6.9}  $[$A special ruled cubic surface$]$ {\rm We consider the surface with equation $$x^2z-xyw+y^3=0.$$
The line $L: x=y=0 $ is a double line for the surface. There are no other singularities. The line $L$ has distinct tangents everywhere except at the point $P:x=y=w=0$, where it has a more complicated singularity. This surface may be regarded as a degenerate case of the ruled cubic surface considered above (see Example \ref{6.6}). It is still a projection of the cubic scroll in $\P^4$, but the conic $\Gamma$ has become two lines.

To investigate the singularity at $P$, we restrict to the affine piece where $z=1$, and the surface $X$ has affine equation $$x^2-xyw+y^3=0\, .$$
Adjoining $u=x/y$ which is an integral element over the affine ring of the surface, we find that the normalization $S$ is just $k[u,w]$, and the map $S \to X$ is given by $x=uy$, $y=u(w-u)$. Lifting the line $L$ to $S$ we find $\Gamma$ is two lines, having equation $u(w-u)=0$, and that our special point $P$ corresponds to the origin $Q:u=w=0$.

The conductor of the integral extension is just $L$, so both $L$ and $\Gamma$ are reduced. We know that $\Cart_Q \Gamma= \Z^2 \oplus k^*$ from Example \ref{5.5}(a), and $\Cart_P L=\Z$. A generator of $\Cart_P L$ is given by $w=0$ on $L$. If we lift $w$ up to $\Gamma$ it intersects each branch in one point giving  $(1,1)$ in $\Z^2$ of $\Cart_Q \Gamma$. Hence $$
\Apic(\Spec \O_{X,P})=\Z\oplus k^*\, .$$

Going back to the projective surfaces $X$ and $S$ we know that $\Pic S=\Z \oplus \Z$, $\Pic \Gamma=\Z \oplus \Z$, $\Pic L=\Z$, so $\Pic X=\Z$ generated by the hyperplane class $H$, whose image in $S$ is $(2,1)$ in the notation of \cite[\S6]{H}. Now by Theorem \ref{TG}, in analogy with the ordinary ruled cubic surface above (cf. \cite[6.3]{H}), an element of $\Apic X$ can be represented by a $4$-tuple $(a,b,\alpha, \lambda)$ where $a,b \in \Z$, $\alpha \in {\rm Div}\, L$, $\lambda \in k^*$ with the condition that ${\rm deg} \, \alpha= a-2b$.
}
\end{Example}
\s

\begin{Example}\label{6.9}$[$The cone over a plane cuspidal curve$]$ {\rm We consider the cubic surface $X$ given by the equation $$y^2z-x^3=0$$ in $k[x,y,z,w],$ which is the cone with vertex $P=(0,0,0,1)$ over a cuspidal curve in the plane $w=0$. The normalization $S$ is obtained by setting $t=\frac{yz}{x}$. It is a cubic surface in $\mathbb P^4$ that is the cone over a twisted cubic curve with vertex $Q$. The singular locus of $X$ is the line $L: x=y=0$. 
Its inverse image in $S$ is the double line $\Gamma$ defined by $x=y=t^2=0$. Since $\Gamma$ is a planar double line, $\Pic L \to \Pic \Gamma$ is an isomorphism. Furthermore $H^0(\O_L^*)=k^*=H^0(\O_{\Gamma}^*)$, $X$ and $S$ are both integral, hence by Proposition \ref{injective}, $\Pic X \to \Pic S$ is an isomorphism as well.

In place of the second sequence in Theorem \ref{TG}, since $S$ is not smooth, we must use the group $G={\rm ker}(\Apic X \to \Apic(\Spec \O_{S, Q}))$ (see Remark \ref{4.4}). Then Theorem \ref{TG} tells us that $G$ is isomorphic to $\Pic S \oplus \Cart \Gamma/\pi^* \Cart L$. Hence we have an exact sequence
$$
0 \to \Pic S \oplus \Cart \Gamma/ \pi^* \Cart L \lto \Apic X \lto \Apic (\Spec \O_{S, Q}) =\Z/3\Z \to 0.
$$
Here $\Pic S=\Z$, generated by three times a ruling, and for each point $x\in L$, $\Cart_x \Gamma/\pi^*\Cart_x L$ is just a group isomorphic to $W=y^{-1}k[y^{-1}]$ since $\Gamma$ is a double line (see Example \ref{5.5}(e)).
}
\end{Example}

\s
\begin{Example}\label{new6.11} $[$A quartic surface with a double line,  \cite[XII, \S 8, p. 471]{P}$]$ {\rm (This is also the surface used by Gruson and Peskine \cite{GP} and \cite{HB}  in their construction of curves in $\P^3$ with all allowable degree and genus.) 

This surface can be obtained by letting $S$ be $\P^2$ with nine points $P_0,  \ldots, P_8$ blown up. We take the points $P_i$ in general position, so that there is a unique cubic elliptic curve $\Gamma$ passing through them. Then $\Pic S=\Z^{10}$ generated by a line from $\P^2$ and the exceptional curves $E_0, \ldots, E_8$. As usual, we denote by $(a;b_0,\ldots,b_8)$ the divisor class $al-\sum b_iE_i$. Thus $\Gamma=(3; 1^9)$. Now take $H=(4;2,1^8)$, i.e. a plane curve of degree 4 with a double point at $P_0$, and passing through $P_1, \ldots, P_8$. Then the complete linear system $|H|$ maps $S$ to a quartic surface $X$ in $\P^3$, whereby 
the curve $\Gamma$ is mapped $2-1$ to a double line $L$ of $X$ with four pinch points. This is a surface $X$ with ordinary singularities and normalization $S$, but it does not arise 
by projection from an embedding of $S$ in some $\P^n$, 
 because the linear system $|H|$ is not very ample on $S$. 

Since $\Gamma$ is an elliptic curve, we have an exact sequence 
$$ 0 \to {\rm Pic}^0 \Gamma \lto \Pic \Gamma \lto \Z \to 0
$$
where ${\rm Pic}^0 \Gamma$ is the Jacobian variety, which in this case is just a copy of the curve itself with its group structure. If we have taken the points $P_i$ in very general position, then the restriction map $\Pic S \to \Pic \Gamma$ will be injective. Dividing by $\pi^* \Pic L =\Z$, generated by the image of the hyperplane class $H$, we see that $\Pic X=\Z$, generated by $H$.  Along the double curve $L$, we have $\Apic (\Spec \O_{X,P})=\Z$ for a general point $P$, or $\Z/2\Z$ at each pinch point.

 }
\end{Example}

\s

\begin{Example}\label{new6.12}$[$A quartic surface with two  disjoint lines as its double locus, \cite[XII, \S 10, type XI, p.490]{P}$]$ {\rm 
The surface has two disjoint double lines as its singular locus $L$, each having four pinch points. The inverse image $\Gamma$ of $L$ in the normalization $S$ will be the disjoint union of two elliptic curves. 

We can construct this surface using an elliptic ruled surface. Following the notations of \cite[V, \S 2]{HAG}, let $\C$ be an elliptic curve. Take $\E=\O_{\C}\oplus \L_0$, where $\L_0$ is an invertible sheaf of degree 0, not isomorphic to $\O_{\C}$, corresponding to a divisor $\e$ on $\C$. We take $S=\P(\E)$. The map $\O_{\C} \to \E$ gives a section $C_0$ of $S$, so that $\Pic S=\Z \cdot C_0 \oplus \Pic \C \cdot f$, where $f$ is a fiber. 

On the surface $S$ we have $C_0^2=0$, $\C_0\cdot f=1$, $f^2=0$. The canonical class is $K_S=-2C_0+\e f$. The surjection $\E \to \O_{\C} \to 0$ defines another section $C_1$ of $S$ that does not meet $C_0$. We have $C_1 \sim C_0 -\e f$, so that $K_S=-C_0 -C_1$.

Now we fix a divisor class $\b$ on $\C$ of degree 2, and take $H=C_0+\b f$ on $S$. Using \cite[V, Ex. 2.11]{HAG} we deduce that the linear system $|H|$ has no base points. Furthermore, $H^2=4$ and $H^0(\O_S(H))=4$, so the linear system $|H|$ determines a morphism $\varphi$ of $S$ to $\P^3$, whose image we call $X$. One can verify (we omit the details) that $H$ collapses the two elliptic curves $C_0$ and $C_1$ to two disjoint lines $L_0$ and $L_1$, and otherwise is an isomorphism of $S-C_0-C_1 \to X-L_0-L_1$. So we denote by $L$ the two disjoint lines $L_0 \cup L_1$ and by $\Gamma$  the two elliptic curves $C_0 \cup C_1$. The $2-1$ coverings $C_0 \to L_0$ and $C_1 \to L_1$ have four ramifications points each., corresponding to four pinch points on each line. In order to find $\Pic X$ we first show that $\Pic S \to \Pic \Gamma$ is injective. Consider a divisor $\eta=nC_0+\a f$ in $\Pic S$. Suppose $\eta$ goes to zero in $\Pic \Gamma$. On $C_1$ the restriction of $C_0$ is zero, since $C_0$ and $C_1$  do not meet.  Now the image of  $\eta$ in $\Pic C_1 \subset \Pic \Gamma$ is $\a$, which is therefore zero. Consider the restriction of $\eta=n C_0$ to $\Pic C_0$.  Since $C_0^2$ is the divisor $\e$, the image of $ \eta$ is $n\e$. If we have chosen $\e$ general, in particular a non torsion element in $\Pic \C$, then $n \e=0$ implies $n=0$. Therefore, when we divide by $\pi^* \Pic L$ only $H$ vanishes. Hence $\Pic X=\Z$ generated by $H$.  As before, $\Apic(\Spec \O_{X,P})$ is $\Z$ at a general point of $L$, and $\Z/2\Z$ at a pinch point.

}
\end{Example}

\begin{Example}\label{quartic}{\rm We now consider a more complicated example. Let $X$ be the quartic surface in $\mathbb P^3$ defined by the polynomial 
$f=x^4-xyw^2+zw^3.$ 

Taking partial derivatives one can verify that the singular locus $L_0$ of $X$ is the line $x=w=0$. Since $f\in (x,w)^3$, the line $L_0$ has multiplicity 3 on $X$. 

To find the normalization $S$ of $X$ we proceed as follows. This method was inspired by a computation in Macaulay 2.  First consider $\frac{f}{w^2}$and let $t=\frac{x^2}{w}$. Then we find $t^2-xy+zw=0,$ so $t$ is integral over the coordinate ring of $X$.  Next we consider $\frac{z^2f}{x^3}$ and let $s=\frac{zw}{x}$. We obtain $s^3-ys^2+xz^2=0,$ thus $s$ is integral over $X$. Now the inverse image of $X$ in the projective space $\mathbb  P^5$ with coordinates $x, y, z, w, s, t$ is the surface $S$ with equations
$$
\begin{cases} tw-x^2 \\ t^2-xy+sx \\ sx-zw \\ tz+s(s-y) \\tx+w(s-y)\\st-xz\end{cases} $$
We recognize this equations as the $2 \times 2$ minors of the $2\times 4$ matrix 
\[ \left( \begin{array}{cccc}
x & s & w & t\\
t & z & x & y-s \\
\end{array} \right).\] 

Thus $S$ is rational scroll of type $(1 ,3)$, namely an embedding of the rational ruled surface $X_2$ (in the notation of \cite[V, \S2]{HAG}). We have $\Pic S =\Z \oplus \Z$, generated by two lines $\C_0$ and $f$ with intersection $\C_0^2=-2$, $\C_0 \cdot f=1$, and $f^2=0$. The hyperplane section $H$ is $\C_0+3f$.

The pullback $\Gamma_0$ of the singular locus $L_0$ is $\C_0+2f$ where $\C_0=(x,w,t,y-s)$ and $f=(x,w,s,t)$. 

To find the conductor of the integral extension, it is enough to look at any affine piece. So let $z=1$ and look at the affine ring $$A=k[x,y,w](x^4-xyw+w^3).$$
The normalization is the affine piece of $S$ defined by $z=1$. One can eliminate variables and find its affine ring is $B=k[s,y],$  and the map $A \lto B$ is defined by $x=s^2(y-s)$, $w=sx$. We claim the conductor of the integral extension is just the ideal $\mathfrak c=(x^2,xw,w^2)$. To see this, observe that $B$ as an $A$-module is generated by $1,s,s^2$. We have only to verify that the elements of $\mathfrak c$ multiply $s$ and $s^2$ into $A$. For instance, $x^2 \cdot s=xw$, $xw \cdot s=w^2$, etc..

We now take the scheme $L$ in $X$ to be defined by the conductor so that we have an exact sequence 
$$ 0 \to \O_{L_0}^2(-1) \lto \O_L \lto  \O_{L_0} \to 0 $$
and hence on units, using Lemma \ref{exp}
$$ 0 \to \O_{L_0}^2(-1) \lto \O_L^* \lto  \O_{L_0}^* \to 0 $$
It follows that $H^0(\O_L^*)=k^*$ and $\Pic L =\Z$. Also, for any point $P\in L$, we find as in Example \ref{5.5}(e)
$$ \Cart_P L \cong \Z \oplus W^2.
$$
We take $\Gamma \subset S$ to be the pullback of $L$. Then $\Gamma=2\C_0+4f$. From the study of curves on ruled surfaces \cite[V, 2.18]{HAG} one knows that $\Gamma$ is linearly equivalent to an irreducible nonsingular curve on $S$. Hence $H^0(\O_{\Gamma})$, which depends only on the linear equivalence class of $\Gamma$, is just $k$. It follows that $H^0(\O_{\Gamma}^*)=k^*$. The canonical class $K_S$ on $S$ is $-2\C_0-4f$ \cite[V, 2.18]{HAG} so from the adjunction  formula one can compute that $p_a(\Gamma)=1$. It follows that $H^1(\O_{\Gamma})=1$. Now we consider the exact sequence 
$$
0 \to \I \lto \O_{\Gamma} \lto \O_{\Gamma_0} \to 0.
$$
Since $p_a(\Gamma_0)=0$, we have $H^0(\O_{\Gamma_0})=1$ and $H^1(\O_{\Gamma_0})=0$. Now from the exact sequence of cohomology we find $H^0(\I)=0$ and $H^1(\I)=1$. Next, we consider the associated sequence of units
$$
0 \to \I \lto \O_{\Gamma}^* \lto \O_{\Gamma_0}^* \to 0.$$
Taking cohomology we obtain
$$
0\to H^1(\I) \lto \Pic \Gamma \lto \Pic \Gamma_0 \to 0.
$$
An analogous argument comparing $\Gamma_0$ to $(\Gamma_0)_{\rm red}$ shows that $\Pic \Gamma_0= \Pic (\Gamma_0)_{\rm red}=\Z \oplus \Z $
(see Example \ref{5.13}(a)). Hence we compute $\Pic \Gamma= \Z \oplus \Z \oplus k^+$.

With the information acquired so far, we can apply Proposition \ref{injective} and conclude that $\Pic X \lto \Pic S$ is injective. Taking for example $\C_0$ and $f$ as a basis for $\Pic S$, it is clear that $\Pic S \lto \Pic \Gamma$ is surjective. Therefore from the sequence of Theorem \ref{TG}(a), we find that $\Pic X=\Z$ and its image in $S$ consists of those curves whose intersection numbers with the two branches of $\Gamma$, that is $\C_0$ and $f$, is the same. These are the curves $a\C_0+bf$ with $b=3a$, i.e. just the multiplies of $H$. So $\Pic X=\Z \cdot H$.  

Finally, we will compute $\Apic (\Spec \O_{X,P})$ for a point $P\in L$. According to Theorem \ref{TL} this is $\Cart_Q \Gamma/ \pi^* \Cart_P L,$ where $Q$ is the point or points of $\Gamma$ lying over $P$. 

The point $P_0$ defined by $y=0$ in $L$ has a single point of $\Gamma$ above it. All other points $P\in L$ have two points lying over them. If $P$ is general point $\not=P_0$, then we have seen that $\Cart_P L=\Z \oplus W^2.$  The two points of $\Gamma$ lying over $P$ are on lines of multiplicity 2 and 4, respectively, so their $\Cart_{Q_i} \Gamma=\Z \oplus W$ and $\Z \oplus W^3$ respectively.  Thus for a general point $P\in L$, $\Apic (\Spec \O_{X,P})=\Z \oplus W^2$. 
At the special point $P_0$ the situation is a bit more complicated. As before, $\Cart_{P_0} L=\Z \oplus W^2$. Using the exact sequence
$$
0 \to I \lto \O_{\Gamma}^* \lto \O_{\Gamma_0}^* \to 0$$
and Proposition \ref{NRed} we obtain an exact sequence
$$0 \to ( I \otimes K)/I \lto  \Cart_{Q_0}\Gamma \lto  \Cart_{Q_0}\Gamma_0 \to 0.
$$
On the other hand, by the same method,
$$ \Cart_{Q_0}\Gamma_0=\Z^2 \oplus k^*\oplus W.
$$
We can also regard $I$ as a $k[y]$-module. It will have rank 3, so that $(I \otimes K)/I\cong W^3$. Now finally, taking the quotient, we obtain
$$\Apic (\Spec \O_{X,P_0})\cong \Z \oplus k^* \oplus W^2
$$ 
}
\end{Example}


\section{The search for set-theoretic complete intersections} 

We say a curve $\C$ in $\P^3$ is a {\it set-theoretic complete intersection} (s.t.c.i. for short) if there exist surfaces $X$ and $Y$ such that $\C=X \cap Y$ as sets. If the surface $X$ containing $\C$ is already specified, we will also say $\C$ is a {\it set-theoretic complete intersection on $X$}. 

In this section we first give some general results, building on the work of Jaffe and Boratynski \cite{J} \cite{J3}, \cite{J2}, \cite{Jaffe}, \cite{B0}, \cite{B}, and \cite{B1}. If $\C$ is a set-theoretic complete intersection on a surface $X$ having ordinary singularities we show that the genus of $\C$ is bounded below as a function of its degree. 

Then we examine some particular surfaces and search for all possible curves that are set-theoretic complete intersections on these surfaces.  

\subsection*{Bounds on degree and genus}

\begin{Proposition}\label{7.1} Let $\C$ be a curve on a surface $X$ in $\P^3$ that meets the singular locus in at most finitely many points. Then $\C$ is a set-theoretic complete intersection on $X$ if and only if $r\C=m H$ in $\Apic X$, for some $r$ , $m \ge 1$, where $H$ is a hyperplane section.
\end{Proposition}

\demo If $\C=X \cap Y$ as sets, where $Y$ is a surface of degree $m$, then the scheme $X \cap Y$ will be a multiple structure on the curve $\C$. Since $\C$ is a Cartier divisor on $X - {\rm Sing} \, X$, this will be $r \C$ for some $r \ge 1$, so $r \C = m H$ in $\Apic X$.  \QED

\s

\begin{Corollary}\label{7.2} With the hypotheses of Proposition \ref{7.1}, assume in addition that $\C$ is smooth, and that $\C$ is a set-theoretic complete intersection on $X$. Then 
at each singular point $P$ of $X$ lying on $\C$, the curve $\C$ gives a non-zero torsion element in the local ring $\Apic(\Spec \O_{X,P})$.
\end{Corollary}
\demo If $\C$ is a s.t.c.i. on $X$, then $r\C \sim m H$ for some $r, m \ge 1$, showing that $r \C$ is a Cartier divisor. Hence the local contribution of $r \C$ at $P$ is zero, so $\C$ is torsion in $\Apic(\Spec \O_{X,P})$.  Being a smooth curve it cannot be locally Cartier at a singular point hence it is non-zero. 
\QED
\s

\begin{Proposition}\label{7.3} If  $\, \C \subset X$ is a smooth curve that gives a non-zero torsion element in   $\Apic(\Spec \O_{X,P})$ for some point $P \in X$, then the normalization $S$ of $X$ can have only one point $Q$ lying over $P$. 
\end{Proposition}
\demo Since $r \C$ is locally Cartier for some $r \ge 1$, it will be defined locally by a single non-zero divisor $f$ in the local ring $\O_{X,P}$. If $S$ has several points $Q_1, \ldots, Q_s$ lying over $P$, then $f$ will define a curve at each $Q_i$. Thus $\C$ will have several branches and cannot be nonsingular at $P$. (this result is deduce by Jaffe \cite[3.3]{J} from a more general result of Huneke.)
\QED

\s

\begin{Corollary}\label{7.4} With the hypotheses of Corollary \ref{7.2}, assume that the surface  $X$ has only ordinary singularities. Then $\C$ can meet the singular locus of $X$ only at pinch points.
\end{Corollary}
\demo Indeed, the normalization $S$ of $X$ has two points lying over a general point of $L$ and three points lying over a triple point. (cf. Examples \ref{three planes} and \ref{6.6}).
\QED
\s

\begin{Proposition}\label{7.5} With the hypotheses of Corollary \ref{7.4}, 
assume that ${\rm char} \, k=0$. Then already $2\C=X \cap Y$ for some surface $Y$. In this case we say that $\C$ is a {\em self-linked} curve.
\end{Proposition}
\demo We have seen that $\C$ can meet the double curve $L$ only at pinch points. At a pinch point the local $\Apic $ group is $\Z/2\Z$ ( see Example \ref{6.6}), so $2\C$ will be locally Cartier there. Thus $2\C \in \Pic X$. Our hypothesis says that  $r\C\sim m H$ for some $r, m \ge 1$. Thus $2\C$ becomes a torsion element of the quotient group $\Pic X/ \Z \cdot H$. However, in case of ${\rm char}\, k=0$, this group is torsion free, by Lemma \ref{7.6} below so we can conclude that $2\C \sim m H$ for some (other) $m$, in other words, by Proposition \ref{7.1}, already $2\C$ is a complete intersection $X \cap Y$ for some $Y$ and $\C$ is self-linked.
\QED

\s

The statement of Lemma \ref{7.6} is well known. We give the main steps of the proof below because of the lack of a precise reference. A proof is given in Jaffe \cite[13.2]{J2}, where he assumes normality but he does not really use it. The statement is also given in \cite[Ex. 20.7]{HDT} as an exercise. Furthermore, a proof  when $X$ as dimension at least 3 can be found in \cite[IV, 3.1]{HLN} (this is the proof that we mimic below). The method of proof is similar to the one employed in the proof of Theorem \ref{4.5}.

\s
\begin{Lemma}\label{7.6} If $X$ is a surface in $\P^3$ over a field $k$ of characteristic zero, then $\Pic X / (\Z \cdot H)$ is a torsion free abelian group. 
\end{Lemma}
\demo In the proof of  \cite[IV, 3.1]{HLN}, at each stage of the thickening, we have $\Pic X_{n+1} \lto \Pic X_n \lto H^2({\I}_{X}^n/{\I}_{X}^{n+1})$, and the class of $H$ comes from $\Pic X_{n+1}$, so its image in the $H^2$ is zero. In characteristic $0$, this $H^2$ group is torsion free, so any class in $\Pic X_n$ whose multiple is in the subgroup generated by $H$  will also go to zero in $H^2$, and hence will lift to $\Pic X_{n+1}$. Continuing in this way, it lifts all the way to 
$\Pic {\P}^3=\Z$, generated by $H$, so it is already a multiple of $H$.
\QED

\s

The following result can be basically recovered from Boratynski's works (\cite{B0} and \cite{B1}). 
\s

\begin{Proposition}\label{7.7} (Boratynski) Let $\C$ be a smooth curve in $\P^3$ that is self-linked, so that $2\C=X\cap Y$ as schemes, where $X$ and $Y$ are surfaces of degrees $m,n$ respectively. Then there is an effective divisor $D$ on $\C$ such that 
$$
\O_{\C} (D)\cong \omega_{\C}^{-2}(2m+n-8)
$$
\end{Proposition}
\demo Let $\C'$ be the scheme $2\C$. Then there is an exact sequence 
$$ 0 \to \L \lto \O_{\C'} \lto \O_{\C} \to 0
$$
where $\L$ is an invertible sheaf on $\C$. One knows from an old theorem of Ferrand that $$\L\cong \omega_{\C} (-m-n+4).$$ This follows also from linkage theory (see for instance \cite[4.1]{H}), since $\L=\I_{\C,\C'}$ which is just ${\mathcal Hom} (\O_{\C}, \O_{\C'})$. But $\omega_{\C}\cong {\mathcal Hom} (\O_{\C}, \omega_{\C'})$ and $\omega_{\C'}\cong \O_{\C'} (m+n-4)$ since $\C'$ is the complete intersection of surfaces of degrees $m$ and $n$. Thus $\L\cong \omega_{\C}(-m-n+4)$. 

On the other hand, $\L$ is a quotient of $\I/\I^2$, where $\I=\I_{\C}$ is the ideal sheaf of $\C$ in $\P^3$. Since $\C$ lies on the surface $X$ of degree $m$ , there is a natural map $\O_{\C}(-m) \to \I/\I^2$, whose image maps to zero in $\L$. Hence there is an effective divisor $D$ on $\C$ such that 
$$
0\to \O_{\C}(-m+D) \lto \I/\I^2 \lto \L\to 0
$$
is exact. 

Now from the exact sequence 
$$0 \to \I/\I^2 \lto \Omega^1_{\P^3|\C} \lto \omega_{\C} \to 0
$$
of \cite[II, 8.17]{HAG}, taking exterior powers, we find 
$$
\wedge^2(\I/\I^2) \cong \omega_{\C}^{-1}(-4).$$
But the above sequence with $\L$ implies that 
$$\wedge^2(\I/\I^2) \cong \L(-m+D).$$
Therefore $\L \cong \omega^{-1}_{\C}(m-4-D).$ Combining with the earlier expression for $\L$ gives the result.
\QED

\s

\begin{Theorem}\label{7.8} Let $\C$ be a smooth self-linked curve, so that $2\C$ is the intersection of surfaces of degrees $m$ and $n$ . Then letting $d={\rm degree} \, \C$ and $g={\rm genus} \,  \C$, we have 
$$ d(m+n-7) \le 4g-4$$ and
$$
g \ge \frac{\sqrt{2}}{2} d^{\frac{3}{2}} -\frac{7}{4} d+1\, .$$
\end{Theorem}
\demo
In the proof of the previous proposition we saw that $\O_{\C}(-m+D)$ is a submodule of $\I/\I^2$. Hence $\O_{\C}(m-D)$ is a quotient of its dual $\N$, the normal bundle of $\C$ in $\P^3$. Since $\C$ is smooth, $\N$ is a quotient of ${\mathcal T}_{\P^3|\C}$, which in turn is a quotient of $\O_{\C}(1)^4$. Hence $\N(-1)$ is a sheaf generated by global sections, and therefore the same holds for $\O_{\C}(m-D)$. We conclude that the degree of this sheaf on $\C$ is non negative, i.e. ${\rm deg} \, D \le dm$. Combining with the expression for $D$ in the statement of Proposition \ref{7.7} we find 
$$
-4g+4+d(2m+n-8)\le dm\, ,$$
which gives the first inequality. 

For the second inequality, we use the fact that $m+n \ge 2\sqrt{mn}$ and $mn=2d$ since $2\C=X\cap Y$. Substituting and solving for $g$ gives the result.
\QED

\s
\begin{Corollary}\label{7.9} If $\C$ is a smooth curve that is a set-theoretic complete intersection on a surface $X$ having at most ordinary singularities and ordinary nodes, and $\C \not\subset {\rm Sing} \, X$, and ${\rm char} \, k=0$, then the inequality of Theorem \ref{7.8} hold, taking $m={\rm deg}\, X$ and $n=\frac{2d}{m}$.
\end{Corollary}

\s

\begin{Corollary}\label{7.10} With the hypotheses of Corollary \ref{7.9} we find $g \ge d-3$ except possibly for the pairs $(d,g)=(8,3), (10,6)$, which we are unable to eliminate. 
\end{Corollary}
\demo For $d\ge 11$ the result follows from the second inequality of Theorem \ref{7.8}. For $d\le 10$ we treat case by case, considering the possible $m,n$ for which $mn=2d$. If one of $m$ or $n$ is equal to 2, then $\C$ lies on a quadric surface, which must be a cone. In this case, if $d$ is even $=2a$, then $g=(a-1)^2$; if $d$ is odd $=2a+1$, then $g=a(a-1)$ \cite[Ch. 3 Ex. 5.6 ]{HAG}. In all cases $g \ge d-3$. If one of $m$ or $n$ is equal to 3, this result has been proved by Jaffe \cite[3.1]{J}. Thus we may assume $m,n\ge4$. This leaves only two cases, $d=8$ and $d=10$, in which cases our bound gives $g\ge 3$ and $g\ge 6$ respectively. In case $d=8$, then $m=n=4$, the curve $\C$ on $S$ is linearly equivalent to $2H$ so $g_{\C}=2g_{H}+3$, which is an odd number. This eliminates $(8,4)$. 
\QED

\s

\begin{Remark}\label{7.11} {\rm Corollary \ref{7.10} was proved by Jaffe for curves on cubic surfaces \cite[3.1]{J}, for surfaces having only ordinary nodes \cite[4.1]{J}, and for curves on cones \cite[5.1]{J}. Boratynski \cite[3.4]{B} proved the special case of genus zero, i.e. a smooth rational curve of degree $d\ge 4$ cannot be a s.t.c.i. on a surface of the type described in Corollary \ref{7.10}, under the additional hypothesis that for $d\ge5$ the rational curve is general in its Hilbert scheme. }
\end{Remark}

\s
\s
\s

\subsection*{Examples and existence} 

\
\s

\bigskip

Now we will study some particular surfaces in $\P^3$, with the intention of finding all curves that are s.t.c.i. on them. We preserve our earlier notation: $X$ will be a surface in $\P^3$, $L$ its singular locus , $S$ its normalization, $\Gamma$ the inverse image of $L$, and we look for smooth curves $\C\subset X$, meeting $L$ in only finitely many points, such that $\C$ is s.t.c.i. on $X$. We denote by $\C'$ the support of the inverse of $\C$ in $S$,
which will be a smooth curve on $S$ isomorphic to $\C$.

\begin{Example}\label{7.14}[The ruled cubic surface (cf. Example \ref{6.6})] {\rm We know that $\C$ can meet $L$ only at pinch points (Corollary \ref{7.4}), so $\C'$ can meet $\Gamma$ only at the two ramification points. In order for $\C$ to be smooth, $\C'$ must meet $\Gamma$ transversally at these points, so the intersection multiplicity $\C'_{\cdot}\Gamma \le 2$. Using the notation of \cite[\S6]{H}, the divisor class of $H$ on $S$ is $(2,1)$ and that of $\Gamma$ is $(1,0)$. Since $r\C'=mH$ for some $r,m$, the class $(a,b)$ of $\C'$ must satisfy $a=2b$. But since  $\C'_{\cdot}\Gamma =a\le 2$, there is only one possibility, namely $\C'=(2,1)=H$ in $\Pic S$. Thus we see that if there is a s.t.c.i. curve $\C$ on $X$, it must be a twisted cubic curve, with $\C'\sim H$ on $S$. To see if such curves $\C$ exist, we ask for a smooth curve $\C'$ in the linear system $|H|$ on $S$ that meets $\Gamma$ at the two specified ramification points. Remembering that $S$ is isomorphic to a plane $\P^2$ blown up at one point $P$, and that $\Gamma$ corresponds to a line $\ell$ not containing $P$, we ask for a conic $\C''$ in the plane containing $P$ and meeting $\ell$ at two specified points. These exist, so we conclude that there are twisted cubic curves $\C$ on $X$ such that $2\C$ is a complete intersection on $X$, and that these are the only smooth s.t.c.i. curves on $X$.

}
\end{Example}
\s

\begin{Example}\label{7.15}[The special ruled cubic surface (cf. Example \ref{new6.9})] {\rm In this case the normalization $S$ is the same cubic scroll in $\P^4$ as in the previous example, but $\Gamma$ is now two lines meeting at the point $Q$ that lies over the special point $P$ on $L$. If $\C$ is a s.t.c.i. curve on $X$, then by Proposition \ref{7.3} it can meet $L$ only at $P$, and since $\C$ is smooth, $\C'$ must meet $\Gamma$ at $Q$ without being tangent to either of the two lines of $\Gamma$. As in the previous example, the class of $H$ in $\Pic S$ is $(2,1)$, the class of $\Gamma$ is $(1,0)$, the intersection number $\C'_{\cdot}\Gamma = 2$, and we find that the only possibility is $\C' \sim H,$ so $\C$ will be a twisted cubic curve. 

There is one well-known example of such a curve \cite[Note, p. 381]{HAJ}, namely (modulo slight change of notation) the twisted cubic curve $\C$ whose affine equation in the open set $z=1$ is given parametrically by 
$$
\begin{cases} x=t^3 \\y=t^2 \\ w=2t \end{cases} .$$
For this curve, $2\C$ is the complete intersection of $X$ with the surface $w^2=4yz$.

If we modify these equations by inserting a parameter $\lambda \not=0, \pm i,$ then the curve $\C_{\lambda}$ defined by 
$$
\begin{cases} x=\lambda t^3 \\y=t^2 \\ w=\frac{\lambda^2+1}{\lambda}t \end{cases} $$
is another smooth twisted cubic curve lying on $X$ and passing though the point $P$, and $\C'_{\lambda}$ on $S$ is still the linear system $|H|$.

To see if $\C_{\lambda}$ is a s.t.c.i. on $X$, we must find out if $\C_{\lambda}$ gives a torsion element in the local group $\Apic (\Spec \O_{X,P})$, which, according to Example \ref{new6.9}, is isomorphic to $\Z\oplus k^*$. Since $\C_{\lambda}$ meets each branch of $\Gamma$ just once, the $\Z$ component of $\C_{\lambda}$ in $\Apic$ is 0. 

To find the element of $k^*$ representing $\C_{\lambda}$ in $\Apic$, we recall from Example \ref{5.5}(a) that if $A \cong k[[x,y]]/(xy)$, then an ideal $(x^r+ay^s)$ gives the element $(r,s;a)$ in $\Cart A\cong \Z^2 \oplus k^*$. Since $\Gamma$ has local equation $u(w-u)=0$, it is convenient to make a change of variables $v=w-u$. then $\Gamma$ is defined by $uv=0$, and the curve $\C'_{\lambda}$ is defined by 
$$
\begin{cases} u=x/y=\lambda t \\v=w-u=\lambda^{-1}t\end{cases},$$
so the local equation of $\C_{\lambda}$ is $u-\lambda^2v$, and we find that $\C_{\lambda}$ gives the element $\mu=-\lambda^2$ in $k^*$. The torsion elements in $k^*$ are just the roots of unity. Note that the first curve we discussed above, which corresponds to $\lambda=1$, gives $\mu=-1$, which is a torsion element of order 2 in $k^*$, confirming that $2\C$ is a complete intersection. 

For the other curves $\C_{\lambda}$ in this algebraic family, we see that for a general $\lambda$, the curve $\C_{\lambda}$ is not a s.t.c.i. on $X$, but if $\lambda$ is a root of unity, then $\C_{\lambda}$ will be a s.t.c.i. For any $n\ge 2$, taking $\mu$ to be a primitive $n^{th}$ root of unity, and solving $\lambda^2=-\mu$, we find a curve $\C_{\lambda}$ that is a s.t.c.i. on $X$ of order $n$, but of no lower order. Note, by the way, that even though the $\C'_{\lambda}$ are all linearly equivalent to $H$ on $S$, the curves $\C_{\lambda_{1}}$ and $\C_{\lambda_{2}}$ on $X$ are not linearly equivalent in $\Apic X$ (unless $\lambda_{1}=-\lambda_{2}$, in which case the curves  $\C_{\lambda_{1}}$ and $\C_{\lambda_{2}}$ are the same, replacing $t$ by $-t$) because this class, in the notation of Example \ref{new6.9}, is $(2,1,0,\mu)$.
}
\end{Example}

\s

\begin{Remark}{\rm In particular the set of curves in an algebraic family on $X$ that are s.t.c.i. need not be either open or closed.}
\end{Remark}
\s

\begin{Example}\label{7.17}[The Steiner surface (cf. Example \ref{Steiner})]   {\rm Here the singular locus $L$ is three lines meeting at a point, and each having two pinch points. The normalization $S$ is isomorphic to $\P^2$, so $\Pic S\cong \Z$, generated by a line $\ell$. The hyperplane section $H$ on $X$ corresponds to the divisor class $2\ell$ on $S$. The curve $\Gamma$ on $S$ corresponds to three lines forming a triangle in $\P^2$.

A s.t.c.i. curve $\C$ in $X$ can meet $L$ only at pinch points (Corollary \ref{7.4}), so the curve $\C' \subset S$ meets $\Gamma$ only at the ramification points. Thus $\C'_{\cdot}\Gamma \le 6$. On the other hand, $\C' \sim a \ell$ for some $a\ge 1$, and $\ell_{\cdot}\Gamma =3$, so $a=1$ or 2. Thus we are looking for a line or a conic in $\P^2$ that meets the triangle $\Gamma$ only in the six ramification points. 

Using the notation of Example \ref{Steiner}, we can compute the ramification points. On the line $t=0$, the image of $(0,1,\lambda)$ is $(0,0,1,\lambda+ \frac{1}{\lambda})$. Thus $\lambda$ and $\frac{1}{\lambda}$ have the same image, and the ramification points are where $\lambda=\frac{1}{\lambda}$, namely, $\lambda=\pm 1$.  Thus the six ramification points in the $(t,u,v)$-plane are $(0,1,\pm1), (1,0,\pm1),$ and $(1, \pm1, 0)$. A line (not equal to $t=0, u=0$ or $v=0$) passing through two of these passes through a third and in this way we obtain four conics $\C$ in $X$ for each of which $2\C$ is a hyperplane section (cf. \cite[p.474]{P}). 

On the other hand, since these six ramification points are aligned 3 at the time, there is no smooth conic passing through all 6. Thus the four conics just mentioned are the only s.t.c.i. curves on $X$.

 }
\end{Example}

\s

\begin{Example}\label{7.18}[A rational quartic surface with a double line (cf. Example \ref{new6.11})]   {\rm In this case $\Gamma$ is an elliptic curve mapping 2-1 to the line $L$, so there are four ramification points. We have $H_{\cdot}\Gamma=2$, so we look for curves $\C'\sim H$ or $\C'\sim 2H$ meeting $\Gamma$ only in the  ramification points.  In the linear system $|H|$ on $S$, every divisor meets $\Gamma$ in a pair of the involution $\sigma$ defining the map $\Gamma \to L$. Thus it cannot meet $\Gamma$ in two ramification points. The divisor $2H$ meets $\Gamma$ in the linear system $2\sigma$. If we compare the map $\Gamma \to L$ to the standard double covering of the $x-$axis by the curve $y^2=x(x-1)(x-\lambda)$ \cite[Chapter IV]{HAG}, then $\sigma$ is just pairs of points that add up to 0 in the group law on the cubic curve, and the four ramification points are the points of order two in the group law. These form a subgroup isomorphic to the Klein four group $\Z/2\Z \oplus \Z/2\Z$, and the sum of all four elements of this group is zero. Hence the sum of the four ramification points is in the linear system defined by $2H$ on $\Gamma$.
Now a straightforward computation of cohomology groups shows that $H^0(\O_S(2H)) \to H^0(\O_{\Gamma}(2H))$ is surjective, so there exist curves $\C' \sim 2H$ meeting $\Gamma$ just at the four ramification points, and it is easy to see, using Bertini's theorem, that we may take $\C'$ to be smooth. Thus we find smooth curves $\C$ of degree 8 and genus 7 on $X$ for which $2\C$ is a complete intersection, and these are the only s.t.c.i. on $X$.

 }
\end{Example}

\s

\begin{Example}\label{7.19}[A quartic surface with two disjoint double lines (cf. Example \ref{new6.12})  ]  {\rm  Here the singular locus $L$ is two lines, and $\Gamma$ is two disjoint elliptic curves, each having four ramification points. Since $H_{\cdot}\Gamma=4$, we look for curves $\C' \sim H$ or $\C'\sim 2H$ meeting $\Gamma$ only in the  ramification points. Any curve $\C' \sim H$ maps to a singular curve in $X$, so we eliminate this case. Consider $\C'\sim 2H$. As in the previous example,  
the sum of the ramification points on $\Gamma$ is in the linear system induced on $\Gamma$ by $2H$. Again, standard calculations of cohomology together with Bertini's theorem show that we can find smooth curves $\C'\sim 2H$ meeting $\Gamma$ only at the ramification points. These give smooth curves $\C \subset X$ of degree 8 and genus 5 for which $2\C$ is a complete intersection, and these are the only s.t.c.i. on $X$. 
}
\end{Example}

\s

\begin{Example}\label{7.20}[The quartic surface of  Example \ref{quartic}]   {\rm We will show that there are no smooth s.t.c.i. curves on this surface. By a reasoning similar to the previous examples we are looking for a smooth curve $\C'\subset S$ meeting $\Gamma$ only at the special point and not tangent to either branch of $\Gamma$. So its local equation will be $y-\lambda s+$ higher terms, for some $\lambda \not=0,1$. For this computation it will be sufficient to use the reduced line $L_0$ and its inverse image $\Gamma_0$ defined by $s^2(y-s)=0$. We will use Proposition \ref{P2}, and show that the image in $\Cart \Gamma_0/\pi^*\Cart L_0$ of the class of $\C'$ in $\Apic(\Spec \O_{X,P_0})$ is not torsion. In Example \ref{quartic}, we found $\Cart_{Q_0} \Gamma_0=\Z^2 \oplus k^8 \oplus W$. Since $\Cart_{P_0} L_0=\Z$, we are looking at the group $\Z \oplus k^*\oplus W$. Now according to Proposition \ref{NRed} and Example \ref{5.5}(e), we write $y-\lambda s +$higher terms as $y(1-\lambda s y^{-1}+ \mbox{higher})$, and find that the contribution of $\C'$ to $W$ is $-\lambda$. The group law in $W$ is addition. It is a $k$-vector space, and in characteristic zero, it is torsion free. Thus $\C'$ gives a non-zero element, and no multiple is zero, so $\C$ cannot be a s.t.c.i.

}
\end{Example}

\s
Recall that the curves that are self-linked in $\P^3$ are s.t.c.i. and conversely on a surface with ordinary singularities in characteristic zero any smooth s.t.c.i. curve not contained in the singular locus  is self-linked. 
\s

\begin{Remark}\label{7.21}{\rm Summarizing what is known so far, we list the possible degree $d \le 8$ and genus $g$ smooth curves $\C$ that are self-linked in $\P^3$ (always assuming ${\rm char} \, k=0$):
$$
\begin{cases} d=3, g=0^{b}\\d=4, g=1^{a}\\d=5, g=2^{b}\\
d=6, g=3^{c},4^{a}\\
d=7, g=6^{b}\\
d=8, g=5^{d},7^{d}, 9^{a}\end{cases}.$$
These are the only possibilities, except maybe for $(d,g)=(8,3)$, which seems unlikely, but we cannot yet exclude. Notes
$$
\begin{cases} a) \quad \mbox{strict complete intersection}\\
b) \quad  2\C \  \mbox{is a complete intersection on a quadratic cone}\\
c) \quad \mbox{ Gallarati \cite{G} finds this one on a cubic surface with four nodes}\\
d) \quad \mbox{Gallarati finds these on quartic surfaces with nodes. The case of $g=5$ lying on a Kummar}\\
\hspace{.7 cm} \mbox{surface was known to Humbert (1883). We find both cases on quartic surfaces with ordinary }\\
\hspace{.7cm}\mbox{(singularities Examples \ref{7.18}, \ref{7.19}).}
\end{cases}.$$

}
\end{Remark}

\s

\begin{Remark}\label{7.22}{\rm The curve of degree 8 and genus 5 is the first known example of a non-arithmetically Cohen-Macaulay curve in $\P^3$ that is a set-theoretic complete intersection in characteristic zero. (In characteristic $p>0$ there are many \cite{HAJ}).

}
\end{Remark}

\s

\begin{Remark}\label{7.23}{\rm using the same techniques, we looked for smooth s.t.c.i. curves on a surface $X$ that is generic projection of a smooth surface $S$ in $\P^n$, for some well-known surfaces $S$. We assume always that $\C$ is not contained in the singular locus of $X$, and that ${\rm char} \, k=0$. We summarize the results here without the computations.
\begin{enumerate}
\item If $S$ is a del Pezzo surface of degree $n$ with $4 \le n \le 9$, there are no s.t.c.i. curves on $X$ except in the case $n=4$, where there are quartic elliptic curves $\C$ with $2\C$ a complete intersection. Of course these curves are already strict complete intersection in $\P^3$. 
\item If $S$ is the quintic elliptic scroll in $\P^4$, there are no s.t.c.i. curves on $X$.
\item Suppose $S$ is a rational scroll $S_{e,n}$ for any $n > e \ge 0$ of degree $d=2n-e$ and having an exceptional curve $\C_0$ with $\C_0^2=-e$. In case $d=3$, we obtain the ruled cubic surface studied above (Example \ref{7.14}). For all $d\ge 4$ there are no s.t.c.i. curves on $X$.
\item Suppose $S$ is the $n$-uple embedding of $\P^2$ in $\P^N$, for $n \ge 2$. In case $n=2$ we obtain the Steiner surface discussed in Example \ref{7.17}. If $n \ge 3$, there are no smooth s.t.c.i. curves on $X$.
\end{enumerate}
Our conclusion is that smooth s..t.c.i. curves on surfaces with ordinary singularities are rather rare!

 }
\end{Remark}

\s

\end{document}